\documentclass[11pt]{article}
\usepackage{amssymb,epsfig,amsmath}
\usepackage{subfigure}
\usepackage{appendix}
\usepackage{rotating}
\newtheorem{lemma}{Lemma}[section]
\newtheorem{theorem}{Theorem}[section]

\newtheorem{coro}{Corollary}[section]

\newtheorem{remark}{Remark}[section]

\numberwithin{equation}{section}

\newcommand{\bsq}{\vrule height .9ex width .8ex depth -.1ex}
\newcommand{\RR}{{\mathbb R}}

\newcommand{\sQ}{{\cal Q}}

\newcommand{\hsp}{\hspace*{\parindent}}
\newcommand{\eeq}{\end{equation}}

\newcommand{\ra}{\rightarrow}
\newcommand{\deq}{\stackrel{\rm d}{=}}
\newcommand{\beql}[1]{\begin{equation}\label{#1}}
\newcommand{\eqn}[1]{(\ref{#1})}
\newcommand{\beq}{\begin{displaymath}}
\newcommand{\eeqno}{\end{displaymath}}

\newcommand{\bC}{{\mathbf C}}

\newcommand{\qandq}{\quad\mbox{and}\quad}

\newcommand{\qforq}{\quad\mbox{for}\quad}

\newcommand{\qasq}{\quad\mbox{as}\quad}
\newcommand{\qinq}{\quad\mbox{in}\quad}

\newlength{\defbaselineskip}
\newcommand{\setlinespacing}[1]%
           {\setlength{\baselineskip}{#1 \defbaselineskip}}
\newcommand{\doublespacing}{\setlength{\baselineskip}%
                           {1.5 \defbaselineskip}}
\newcommand{\singlespacing}{\setlength{\baselineskip}{\defbaselineskip}}

\begin{document}

\doublespacing

\begin{center}

{\large\bf TWO-PARAMETER HEAVY-TRAFFIC LIMITS} \\
{\large\bf FOR INFINITE-SERVER QUEUES} \\

\vspace{0.1in}

 Guodong Pang and Ward Whitt \\
~ \\

 Department of Industrial Engineering and Operations Research \\
 Columbia University \\
 New York, NY 10027-6699 \\
 Email: {\{gp2224, ww2040\}}@columbia.edu \\
\vspace{0.2in}

July 9, 2010 \\
\vspace{0.1in}

\end{center}

\begin{center}
{\em Abstract}
\end{center}

In order to obtain Markov heavy-traffic approximations for
infinite-server queues with general non-exponential service-time
distributions and general arrival processes, possibly with
time-varying arrival rates, we establish heavy-traffic limits for
two-parameter stochastic processes.  We consider the random
variables $Q^e (t,y)$ and $Q^r (t,y)$ representing the number of
customers in the system at time $t$ that have elapsed service
times less than or equal to time $y$, or residual service times
strictly greater than $y$. We also consider $W^r (t,y)$
representing the total amount of work in service time remaining to
be done at time $t+y$ for customers in the system at time $t$.
 The two-parameter stochastic-process limits in the space $D([0,\infty),D)$ of $D$-valued functions in $D$
draw on, and extend, previous heavy-traffic limits by Glynn and
Whitt (1991), where the case of discrete service-time
distributions was treated, and Krichagina and Puhalskii (1997),
where it was shown that the variability of service times is
captured by the Kiefer process with second argument set equal to
the service-time c.d.f.

\vspace{0.3in} \noindent {\em Keywords}:   infinite-server queues,
heavy-traffic limits for queues, Markov approximations,
two-parameter processes, measure-valued processes, time-varying
arrivals, martingales, functional central limit theorems,
invariance principles, Kiefer process.

\section{Introduction}\label{secIntro}
\hsp
One reason
heavy-traffic limits for queueing systems are useful is that they show that non-Markov stochastic processes
describing system performance
can be approximated by Markov stochastic processes under heavy loads.  For a Markov process,
it suffices to know the present state of that stochastic process in order
 to determine the distribution of the stochastic process at future times;
we need no additional information from the past.
With Markov approximations, that remains true approximately.
In applications, the Markov property shows that the proper state has been identified and shows
what needs to be measured in order to understand system performance.

The classic example is the conventional heavy-traffic limit for
the $GI/GI/s$ queue, having $s$ servers, unlimited waiting room,
and independent and identically distributed (i.i.d.) service times
independent of a renewal arrival process. The standard description
of system state is the number of customers in the system at time
$t$, which we will call the queue length and denote by $Q(t)$.
With non-exponential interarrival and service times, the
stochastic process $\{Q(t): t \ge 0\}$ is not Markov.  Then the
future evolution at any time depends on the elapsed interarrival
time and the elapsed service times of all customers being served.
However, the conventional heavy-traffic limit, in which the
traffic intensity approaches the critical value $1$ from below
while the number of servers remains fixed, shows that the
queue-length process $\{Q(t): t \ge 0\}$ is approximately equal to
a Markov process, in particular, reflected Brownian motion, under
heavy loads \cite{IW70a, IW70b, W02}. In fact, the interarrival
times and service times need not come from independent sequences
of i.i.d. random variables. Instead, it suffices to have the
associated partial sums, or equivalently, the associated counting
processes satisfy a FCLT. Moreover, the Markov property of the
limit extends to conventional heavy-traffic limits for networks of
queues \cite{HR81}.

The situation is very different for many-server heavy-traffic
limits when the service-time distribution is non-exponential,
either with $s = \infty$ or $s \ra \infty$.  In this paper, we
will consider the case in which $s = \infty$, i.e., the
$G/GI/\infty$ model with i.i.d. service times independent of a
general arrival process, where heavy traffic is achieved by
letting $\lambda \ra \infty$, while the service-time distribution
is held fixed. However, the problem is relevant more generally
with many servers, where $s \ra \infty$ as $\lambda \ra \infty$
with $s - \lambda = O(\sqrt{\lambda})$, as in \cite{HW81}.
For infinite-server models, we index the stochastic
processes by the arrival rate $\lambda$. We are interested in the
infinite-server model both for its own sake and as an
approximation for many-server queues.  In fact, heavy-traffic
limits for infinite-server models can play a role in
characterizing the heavy-traffic limits for corresponding
many-server models, as shown by Reed \cite{R07a, R07b} and \cite{MM09, PR08}.

With infinitely many exponential servers, we again obtain
Markov diffusion limits, as first shown by Iglehart \cite{I65} for the $M/M/\infty$ model;
see \cite{PTW07} for a review.
A systematic way to extend the limit to general arrival processes is given in \S7.3 of \cite{PTW07}.
However, with non-exponential service times,
the established heavy-traffic limit for $Q(t)$ is
{\em not} Markov.  As first shown by Borovkov \cite{B67}, and further discussed in \cite{I73, W82, L88, GW91, KP97},
the limit process is Gaussian, which implies that the distribution of $Q(t)$ itself is approximately
normal, but the limiting Gaussian stochastic process is non-Markov,
unless the service times are exponential (plus a minor additional case, \cite{G82, KP97}).

We consider a stochastic process characterizing the system state for which
the associated heavy-traffic limit process is Markov.
We do so in two ways:  First, we
consider the two-parameter stochastic process $\{Q^e(t,y): t \ge 0, y \ge 0\}$, where
$Q^e(t,y)$ represents the number of customers in the system at time $t$ with elapsed service times less
than or equal to time $y$.  We do not pay attention to specific customers or servers
but only count the total numbers.
The random quantity $Q^e(t,y)$ is an {\em observable} quantity given the system history up to time $t$.
We recommend that the stochastic process $\{Q^e(t,y): t \ge 0, y \ge 0\}$ be used in models and measured in practice.
Ways to exploit such ages for control are discussed in \cite{DW97}.

So far, we have used elapsed service times, because they are directly observable.  We can equally well work
with residual service times, and consider the process $Q^r (t,y)$ counting the number of customers in the system at time $t$
with residual service times strictly greater than $y$.
  With i.i.d. service times having c.d.f. $F$, we can go from one formulation to the other.
If the elapsed service time is $y$, then the residual service time has distribution $F_y (x) \equiv F(x+y)/F^c(y)$ for $x \ge 0$,
where $F^c (y) \equiv 1 - F(y)$.
If the service times are learned when service begins, then both $Q^r (t,y)$ and $Q^e(t,y)$ are directly observable.
Otherwise, elapsed service times correspond to what we observe, while residual service times
represent the future load, whose distribution we may want to describe.

We regard $\{Q^e (t,\cdot): t \ge 0\}$ and $\{Q^r (t,\cdot): t \ge 0\}$ as function-valued stochastic processes, in particular, random elements of
the function space $D_D$; see \S \ref{secD}.
Since the functions $Q^e(t,y)$ ($Q^r(t,y)$) are
nondecreasing (nonincreasing) in $y$, we can also regard $Q^e(t, \cdot)$ and $Q^r(t, \cdot)$ as measure-valued processes,
but we will work in the framework $D_D$.

For the $M/GI/\infty$ model, it is easy to see that the stochastic process $\{Q^e(t,\cdot): t \ge 0\}$
is a Markov process;
\cite{EMW93, MW93} (where references to earlier work are given).
The key idea, expressed in the proof of Theorem 1 of \cite{EMW93},
is a Poisson-random-measure representation.
For the more general $GI/GI/\infty$ model, having a non-Poisson renewal arrival process,
the stochastic process $\{Q^e(t,\cdot): t \ge 0\}$ is in general not Markov from that perspective,
because the future evolution also depends on the elapsed interarrival time.
The Markov property is violated more severely when the arrival process is not renewal.
However, just as for the $G/GI/s$ model discussed above,
the heavy-traffic limit for the arrival process typically does have independent increments,
so this non-Markovian aspect disappears in the heavy-traffic limit.
In the limit, $Q^e (t,y)$ for the $G/GI/\infty$ model
is asymptotically equivalent to what it would be in the corresponding
$M/GI/\infty$ model, except for a constant factor $c_a^2$ to account for the different variance;
see Theorem \ref{thGaussian} and Corollaries \ref{coro1Var} and \ref{coro2Var}.

\paragraph{Proof Strategy.}
Our proof builds on previous work by Glynn and Whitt \cite{GW91} and
 Krichagina and Puhalskii \cite{KP97}.
First, a restricted form of the desired two-parameter
stochastic-process limit was already established in Theorem 3
of \cite{GW91} for the case of service-time
distributions with finite support.  That result is only
stated in $D$ for arbitrary fixed second parameter $y$, but it can be
extended quite easily to the function space $D_D$.
Since distributions with finite support are dense
in the space of all probability distributions, one might consider the matter settled.
However, much
depends on the precise assumptions made about the service-time distribution.  The goal should be to treat general service-time distributions
without any extra conditions.  We should not need to assume that any moments are finite or
that the c.d.f. is continuous or absolutely continuous.

One important feature of \cite{KP97} is that they treat completely general service-time distributions.
However, they do not state limits for
two-parameter queueing processes.  It might seem that it should be a routine extension to do so,
but we show that is {\em not} so,
because a candidate limit process is not a random element of the space $D_D$ for
discontinuous service-time c.d.f.'s, as we explain in Remark \ref{rm3}.
Fortunately, however, the argument in \cite{KP97} can be extended to the two-parameter case
if we restrict attention to continuous service-time c.d.f.'s, which we do.

A key idea in \cite{GW91} is to treat service-time distributions with finite support by
representing them as finite mixtures of
deterministic service-time c.d.f.'s, and then split the arrival
process into corresponding arrival processes associated with
each deterministic service time; see \S 3 of \cite{GW91},
especially, Proposition 3.1.  That step relies on the FCLT for
split counting processes, as in \S 9.5 of \cite{W02}.   The
mixture argument extends quite directly to treat arbitrary
discrete distributions. It also extends to arbitrary
distributions if we can treat continuous service-time c.d.f.'s,
but the proof in
\cite{GW91} does not seem to extend naturally to continuous
service-time c.d.f.'s.

Hence, for the final case of a continuous service-time c.d.f.,
we draw heavily on \cite{KP97}.  Our limits for
continuous service-time c.d.f.'s are extensions of theirs,
obtained using the same function space and the same martingale arguments. The proof
in \cite{KP97} already took a
two-parameter approach and, following Louchard \cite{L88}, showed that it is fruitful to view
the service times through the associated sequential empirical
process (in \eqn{E:Kn} below).  They showed that a scaled version of
the sequential empirical process converges to the two-parameter
standard Kiefer process, with the service time c.d.f. in the
second argument (see \eqn{b6} below).  This convergence was
established in the space
$D_D$; see \S \ref{secD}.

\paragraph{Other Related Literature.}  As noted in \cite{KP97},
the relevance of the two-parameter Kiefer process for the infinite-server queue was first observed by Louchard \cite{L88}.
 The results here were briefly outlined in
\S6.4 in our survey \cite{PTW07}.  (The first drafts of this paper were written at that time.)
Related fluid limits for measure-valued processes
have since been obtained in \cite{KR07, KR08, Z09}.
However, the first fluid limit for two-parameter processes for this model evidently was the fluid limit in \S 6 of \cite{WW06}
for the discrete-time version of that more general $G_{t}(n)/GI/s+GI$ model, having both time-dependent and state-dependent arrivals.
Decreusefond and Moyal \cite{DM08} established a FCLT for the $M/GI/\infty$ model.
In contemporaneous work, Reed and Talreja \cite{RT09} extend the
result in \cite{DM08} to the $G/GI/\infty$ model and show that the limit process $\hat{Q}^e$ can be
regarded as an infinite-dimensional (distribution-valued) OU process,
thus proving that the limit process $\{\hat{Q}^e (t , \cdot): t \ge 0\}$ is a Markov process.
In these other papers, like \cite{GW91}, there are extra regularity conditions on the service-time c.d.f.
Moreover, the alternative spaces admit fewer continuous functions.

\paragraph{Organization of this paper.}
We start with preliminaries in \S \ref{secPrelim}.
In \S \ref{secMain} we state our main results, focusing only on new arrivals (ignoring any customers initially in the system).
In \S \ref{secChar} we characterize the limit processes.
In \S \ref{secInitial} we treat the initial conditions, and treat all customers in the system.
In \S \ref{secProof} we prove the main theorem:  Theorem \ref{QFCLT}, focusing on the case of continuous service-time distribution.
In \S \ref{secContinuity} we prove the continuity of the representation of some key processes in the space $D_D$.
In \S \ref{secTight} we continue the proof by establishing tightness of the key processes.
In \S \ref{secFiDi} we complete the proof by establishing convergence of the finite-dimensional distributions.
There is also a longer version of the present paper \cite{PW09} available from the authors' web sites.  It has a longer introduction;
it shows how known results for
the special case of exponential service times can be derived from our formulas;
it presents supporting technical details, including basic facts about the
Brownian sheet, the Kiefer process, two-parameter stochastic integrals, tightness criteria in the space $D_D$ and some detailed calculations.

\section{Preliminaries}\label{secPrelim}

\subsection{Initial Conditions and Assumptions}

It is convenient to treat the congestion experienced by customers initially in the system separately
from the congestion experienced by new arrivals, because they usually can be regarded as being asymptotically independent.
Thus we first focus only on new arrivals and then later treat the initial conditions in \S \ref{secInitial}.

\paragraph{Assumptions for the Arrival Processes.}

We consider a sequence of $G/GI/ \infty$ queues indexed by $n$, where the arrival rate is increasing in $n$.
For the $n^{\rm th}$ system, let
$A_n(t)$ be
the number of arrivals by time $t$ and $\tau^n_i$ the time of the $i^{\rm th}$ arrival.

We assume that the sequence of arrival processes satisfy a FCLT, specified below.
All single-parameter continuous-time stochastic processes are assumed to be random elements of the
function space $D \equiv D([0,\infty), \RR)$ with the usual Skorohod $J_1$ topology \cite{B68, W02}.
 Convergence $x_n \ra x$ as $n \ra \infty$ in the $J_1$ topology is equivalent to uniform convergence
on compact subsets (u.o.c.) when the limit function $x$ is continuous.
Throughout, we will have a bar, as in $\bar{A}_n (t)$, to denote the law of large number (LLN) scaling (as in \eqn{b1} below)
and a hat, as in $\hat{A}_n (t)$, to denote the central limit theorem (CLT) scaling (as in \eqn{b2} below).

\paragraph{Assumption 1:  FCLT.}
There exist:  (i) a \textit{continuous} nondecreasing deterministic real-valued function $\bar{a}$ on $[0,\infty)$ with $\bar{a}(0) = 0$
and (ii) a stochastic process $\hat{A}$ in $D$ with continuous sample paths, such that
\beql{b2}
\hat{A}_n (t) \equiv n^{-1/2} (A_n (t) - n \bar{a}(t)) \Rightarrow \hat{A} (t) \qinq D \qasq n \ra \infty.~~~\bsq
\eeq

As an immediate consequence of Assumption 1, we have an associated functional weak law of large numbers (FWLLN)
\beql{b1}
\bar{A}_n (t) \equiv \frac{A_n (t)}{n} \Rightarrow \bar{a}(t) \qinq D \qasq n \ra \infty.
\eeq
In order to obtain a limiting Markov process we will also assume that the limiting stochastic process $\hat{A}$ has independent increments,
but we will obtain limits more generally.

\paragraph{The Standard Case.}
The standard case in Assumption 1 has special $\bar{a}$ and $\hat{A}$.
For the FWLLN limit, the standard case is $\bar{a}(t) = \lambda t, t \ge 0$ for some positive constant $\lambda$, which
corresponds to an arrival rate of $\lambda_n \equiv \lambda n$ in the $n^{\rm th}$ system,
but our more general form allows for time-varying arrival rates as in \cite{EMW93, MW93, MMR98}.

For the FCLT limit $\hat{A}$, the standard case is BM.  That occurs when the arrival processes
are scaled versions of a common renewal process with interarrival times having mean $\lambda^{-1}$ and
SCV $c_a^2$.  Then
$\hat{A} (t) = \sqrt{\lambda c_a^2} B_a (t)$, where $B_a$ is a standard BM.  Of course, the convergence to BM in \eqn{b2}
holds much more generally, e.g., see Chapter 4 of \cite{W02}.  Except for the SCV $c^2_a$, in the standard case Assumption 1
makes the arrival process asymptotically equivalent to a Poisson process.
Thus, in the standard case, the limiting results will be identical to the limit for the $M/GI/\infty$ model when $c_a^2 = 1$,
and very similar for $c_a^2 \not= 1$.
Actually, there is an important structural difference when $c_a^2 \not= 1$, which we discuss in \S \ref{secChar}.

\paragraph{Assumptions for the Service Times and the Empirical Process.}

\paragraph{Assumption 2:  a sequence of i.i.d. random variables.}
We assume that the service times of new arrivals come from a sequence of i.i.d.
nonnegative random variables $\{\eta_i: i \ge 1\}$ with a \textit{general}
c.d.f. $F$, independent of $n$ and
the arrival processes. ~~~\bsq

As in \cite{KP97}, it is significant that our queue-length heavy-traffic limits over finite time intervals
do not require more assumptions about the service-time c.d.f. $F$. It need not have a finite mean.
However, for subsequent results we will need to assume in addition that $F$ has a finite mean $\mu^{-1}$ and even
a finite second moment with SCV $c^{2}_s$.

Krichagina and Puhalskii \cite{KP97} observed that
it is fruitful to view the service times through the two-parameter {\em sequential empirical process}
\beql{E:Kn}
\bar{K}_n (t,x) \equiv \frac{1}{n}\sum\limits_{i = 1}^{\lfloor n t \rfloor}\mathbf{1}(\eta_i \leq x), \quad t \geq 0, \quad x\geq 0,
\eeq
which is directly expressed in the LLN scaling. Here $\mathbf{1}(A)$ is the indicator function.
Since the service times are i.i.d. (without any imposed moment conditions),
we have a FWLLN for $\bar{K}_n$ itself and a FCLT for the scaled process
\beql{b4}
\hat{K}_n (t,x) \equiv \sqrt{n}(\bar{K}_n (t,x) - E[\bar{K}_n (t,x)]) = \frac{1}{\sqrt{n}}\sum\limits_{i = 1}^{\lfloor n t \rfloor}
(\mathbf{1}(\eta_i \leq x) - F(x)),
\eeq
for $t\geq 0$ and $x\geq 0$.

These stochastic-process limits are based on corresponding limits in the case of random variables uniformly distributed on $[0,1]$.
Let $\hat{U}_n (t,x)$ denote the stochastic process $\hat{K}_n (t,x)$ when $\eta_i$ is uniformly distributed on $[0,1]$,
so that $F(x) = x$, $0 \le x \le 1$.
Extending previous results by Bickel and Wichura \cite{BW72}, Krichagina and Puhalskii \cite{KP97} showed that
\beql{b5}
\hat{U}_n (t,x) \Rightarrow U (t,x)  \qinq D([0,\infty), D([0,1], \RR)) \qasq n \ra \infty,
\eeq
where $U(t,x)$ is the {\em standard Kiefer process}; see Cs\"{o}rg\"{o} and R\'{e}v\'{e}sz \cite{CR81},
Gaenssler and Stute \cite{GS79}, and van der Vaart and Wellner \cite{VW96}.
In particular, $U(t,x) = W(t,x) - xW(t,1)$, where $W(t,x)$ is a two-parameter BM (Brownian sheet), so that
$U(\cdot, x)$ is a BM for each fixed $x$, while $U(t, \cdot)$ is a Brownian bridge for each fixed $t$.
The Brownian bridge $B^0$
can be defined in terms of a standard BM $B$ by $B^0 (t) \equiv B(t) - t B(1)$, $0 \le t \le 1$;
it corresponds to BM conditional on having $B(1) = 0$.

It is significant that $\hat{K}_n$ can be expressed as a simple composition of $\hat{U}_n$ with the c.d.f. $F$
in the second component.  We thus have
\beql{b6}
\hat{K}_n (t,x) = \hat{U}_n (t,F(x)) \Rightarrow  \hat{K}(t,x) \equiv U (t,F(x))  \qinq D([0,\infty), D([0,\infty), \RR)),
\eeq
as $n\ra\infty$
without imposing any conditions upon $F$, because $F$ is not dependent on $n$.  Moreover,
the convergence is with respect to a stronger topology on $D_D \equiv D([0,\infty), D([0,\infty), \RR))$;
convergence is uniform over sets of the form $[0,T] \times [0, \infty)$; we have uniformity over $[0,\infty)$ in the second argument.
That will turn out to be important when we treat the remaining-workload process.
As a consequence of the FCLT in \eqn{b6}, we immediately obtain the associated FWLLN
\beql{b7}
\bar{K}_n (t,x) \Rightarrow \bar{k}(t,x) \equiv t F(x) \qinq D_D \qasq n\ra\infty,
\eeq
where again there is uniformity in $x$ over $[0, \infty)$.

\subsection{Prelimit Processes}

Let $Q_n^e (t,y)$ represent the number of customers in the $n^{\rm th}$ queueing system at time $t$
that have {\em elapsed} service times less than or equal to $y$;
let $Q_n^r (t,y)$ represent the corresponding number
that have {\em residual} service times strictly greater than $y$.
Let $Q_n^t (t)$ represent the {\em total number} (the superscript $t$) of customers in the $n^{\rm th}$ queueing system at time $t$.
Clearly, $Q_n^t (t) = Q_n^e (t,t) = Q_n^r (t,0)$, and
\begin{eqnarray}\label{b3}
Q_n^r (t,y) &=& Q_n^e(t+y,t+y) - Q_n^e(t+y,y) = Q_n^{t}(t+y) - Q_n^e(t+y,y), \nonumber \\
Q_n^e (t,y) &=& Q_n^r(t,0) - Q_n^r(t-y,y) = Q_n^t(t) - Q_n^r(t-y,y).
\end{eqnarray}
From \eqn{b3}, it is evident that we can construct all three processes $Q_n^e$, $Q_n^r$ and $Q_n^t$ from
either $Q_n^e$ or $Q_n^r$.
 Observe that  $Q_n^r$ and  $Q_n^e$ can be expressed as
\begin{eqnarray}\label{E:Qn}
Q^r_n (t,y) &=&  \sum\limits_{i = 1}^{A_n(t)} \mathbf{1}( \tau^n_i+ \eta_i> t+y),  \quad t\geq 0, \quad y \geq 0, \\
Q^e_n (t,y) &=&  \sum\limits_{i = A_n(t-y)}^{A_n(t)} \mathbf{1}( \tau^n_i+ \eta_i> t), \nonumber \quad t\geq 0, \quad  0 \leq y \leq t.
\end{eqnarray}
From \eqn{E:Qn}, we see the connection to the sequential empirical process $\bar{K}_n$ in \eqn{E:Kn}.
Indeed, the key observation (following \cite{KP97}) is that we can rewrite the random sums in \eqn{E:Qn} as integrals
with respect to the random field $\bar{K}_n$ by
\begin{eqnarray}\label{E:QnK}
Q^r_n (t,y) &=&
n \int_0^t \int_0^{\infty} \mathbf{1}(s+x > t+y) d \bar{K}_n(\bar{A}_n(s),x), \quad t, y\geq 0,  \\
Q^e_n (t,y) &=&
n \int_{t-y}^t \int_0^{t} \mathbf{1}(s+x > t) d \bar{K}_n(\bar{A}_n(s),x), \quad t\geq 0, \quad 0 \leq y \leq t, \nonumber
\end{eqnarray}
for $\bar{K}_n$ in \eqn{E:Kn}.  These two-dimensional integrals in \eqn{E:QnK} are two-dimensional
Stieltjes integrals.  In the present context, the integrals in \eqn{E:QnK} are understood to be (defined as)
the random sums in \eqn{E:Qn}.

\begin{lemma}{\em $($representation of $Q_n^r$ and $Q_n^e)$}\label{repQn}
The processes $Q^r_n$ and $Q_n^e$ defined in {\em \eqn{E:Qn}} and {\em \eqn{E:QnK}} can be represented as
\beql{E:QnM}
Q^r_n (t,y) = n \int_0^t F^c (t+y-s) \, d\bar{a}(s)
 + \sqrt{n}(\hat{X}^r_{n,1}(t,y) + \hat{X}^r_{n,2}(t,y)),\quad t, y \geq 0,
\eeq
\beql{E:QneM}
Q^e_n (t,y) = n \int_{t-y}^t F^c (t-s) \, d\bar{a}(s)
 + \sqrt{n}(\hat{X}^e_{n,1}(t,y) + \hat{X}^e_{n,2}(t,y)),\quad t \geq 0, \quad 0\leq y \leq t,
\eeq
where
\begin{eqnarray}\label{E:Mn1}
\hat{X}^r_{n,1}(t,y) \equiv  \int_0^t F^c (t+y-s) \, d \hat{A}_n(s),  \quad \hat{X}^e_{n,1}(t,y) \equiv \int_{t-y}^t F^c (t-s) \, d \hat{A}_n(s),
\end{eqnarray}
\begin{equation}\label{E:Mn2}
\hat{X}^r_{n,2}(t,y) \equiv \int_0^t \int_0^{\infty} \mathbf{1}( s+x > t+y)\, d\hat{R}_n(s,x)
        =  - \int_0^t \int_0^{\infty} \mathbf{1}( s+x \leq t+y)\, d\hat{R}_n(s,x),
\end{equation}
\begin{equation}\label{E:Mn2e}
\hat{X}^e_{n,2}(t,y) \equiv \int_{t-y}^t \int_0^{t} \mathbf{1}( s+x > t)\, d\hat{R}_n(s,x)
=  - \int_{t-y}^t \int_{0}^{t} \mathbf{1}( s+x \leq t)\, d\hat{R}_n(s,x),
\end{equation}
\begin{eqnarray} \label{E:Rn}
\hat{R}_n(t,y)&  \equiv & \hat{K}_n (\bar{A}_n (t),y)  =  \frac{1}{\sqrt{n}}\sum\limits_{i =1}^{A_n(t)}(\mathbf{1}(\eta_i \leq y)-F(y)) \\
& = &\sqrt{n} \bar{K}_n (\bar{A}_n (t),y) - \hat{A}_n (t) F(y) - \sqrt{n} \bar{a}(t) F(y), \nonumber
\end{eqnarray}
with the integrals in {\em \eqn{E:Mn1}}, {\em \eqn{E:Mn2}} and {\em \eqn{E:Mn2e}} all defined as Stieltjes integrals
for functions of bounded variation as integrators.
\end{lemma}

\paragraph{Proof.}
Apply \eqn{b4} to get the first relation in \eqn{E:Rn}.  (Right away, from \eqn{b6},
we see that $\hat{R}_n (t, x) \Rightarrow \hat{K} (\bar{a} (t), x)$.)
  Use \eqn{b4} and \eqn{E:Kn} to get the rest of \eqn{E:Rn} and
\begin{eqnarray}\label{E:KnA}
\bar{K}_n(\bar{A}_n(t), x) &=& \frac{1}{n}\sum\limits_{i =1}^{A_n(t)}\mathbf{1}(\eta_i \leq x) \nonumber \\
&=& \frac{1}{\sqrt{n}}\Big[\frac{1}{\sqrt{n}}\sum\limits_{i =1}^{A_n(t)}(\mathbf{1}(\eta_i \leq x)-F(x))\Big] +\frac{1}{\sqrt{n}}\sqrt{n}(\bar{A}_n(t)-\bar{a}(t))F(x)  +  \bar{a}(t)F(x) \nonumber \\
&=&  \frac{1}{\sqrt{n}}\hat{R}_n(t,x) + \frac{1}{\sqrt{n}} \hat{A}_n(t)F(x) +  \bar{a}(t)F(x).
\end{eqnarray}
Combine \eqn{E:QnK} and \eqn{E:KnA} to get \eqn{E:QnM}.
The alternative representation for $\hat{X}_{n,2}(t,y)$ holds because $\hat{K}_n (t, \infty) = 0$ and thus
$\hat{R}_n (t, \infty) = 0$ for all $t$. ~~~\bsq

We will also consider several related processes.
Let $F^e_n (t,\cdot)$ and $F^r_n (t,\cdot)$  represent the {\em empirical age distribution}  and the  {\em empirical residual distribution} at time $t$ in the $n^{\rm th}$ system, respectively, i.e.,
\beql{age}
F^e_n (t,y) \equiv Q^e_n (t,y)/Q^t_n (t), \quad t\geq 0, \quad 0 \leq y \leq t,
\eeq
and
\beql{residual}
F^{r,c}_n (t,y) \equiv 1 -  F^{r}_n (t,y)  \equiv  Q^r_n (t,y)/Q^t_n (t), \quad t\geq 0, \quad y \geq 0.
\eeq
For each $n$ and $t$, $F^e_n (t,\cdot)$ and $F^r_n (t,\cdot)$  are proper  c.d.f.'s.
Let $D_n (t)$ count the number of departures in the interval $[0,t]$; clearly, $D_n (t) \equiv A_n (t) - Q^t_n (t)$ for $t \ge 0$.

We will also consider several processes characterizing the workload in total service time.
For these limits, we will assume that we are in the standard case for the arrival process and impose extra moment conditions
on the service-time c.d.f. $F$.
The total input of work over $[0,t]$ is
\beql{b8}
I_n (t) \equiv \sum_{i=1}^{A_n (t)} \eta_i, \quad t \ge 0.
\eeq
The amount of the workload to have arrived by time $t$ that will be remaining after time $t + y$ is
\beql{E:Vn}
W^r_n (t,y) \equiv \int_{y}^{\infty} Q^r_n (t, x) \, dx, \quad t\geq 0, \quad y \geq 0.
\eeq
Then the total (remaining) workload at time $t$ is $W^t_n (t) \equiv W^r_n (t,0)$.  Finally, the total amount of completed
service work by time $t$ is $C_n (t) \equiv I_n (t) - W_n^t (t)$.

\subsection{The Space $D_D$}\label{secD}

Our limits for two-parameter processes will be in the space $D_D$, which we regard as
a subset of $D([0,\infty), D([0,\infty), \RR))$, where $D \equiv D([0,\infty), S)$, for a separable metric space $S$, is the space of all right-continuous $S$-valued functions with left limits in $(0,\infty)$;
see \cite{B68, W02} for background.
We will be considering the subset of functions $x(t,y)$ which have finite limits as the second argument $y \ra \infty$.
For example, we have $Q^e_n (t,y) = Q^e_n (t,t)$ for all $y > t$ and $Q^r_n (t,y) \ra 0$ as $y \ra \infty$.
We will be using the standard Skorohod \cite{S56} $J_1$ topology on all $D$ spaces, but since all limit
processes will have continuous sample paths, convergence in our space $D_D$
is equivalent to uniform convergence over subsets of the form
$[0,T]\times [0, \infty)$.  (We already observed that we have such stronger uniform convergence over that non-compact set for
$\hat{K}_n$ to the Kiefer process in \eqn{b5}.)
We refer to \cite{TW08} for the convergence preservation of various functions in $D_D$.

For two-parameter processes, one might consider using generalizations of the spaces of two-parameter real-valued functions
considered by Straf \cite{S71} and Neuhaus \cite{N71}, but those spaces require limits to exist at each
point in the domain (subset of $\RR^2$) through all paths lying in each of the four quadrants centered at that point.
That works fine for the sequential empirical process $K_n$, but {\em not} for $Q^r_n (t,y)$.
For example, suppose that the first two arrivals occur at times $1$ and $3$, and that the arrival at time $1$ has
a service time of $2$.  Then limits do not exist along all paths in the southeast and northwest quadrants
at the point $(t,y) = (2,1)$, because there are discontinuities along a negative $45^{o}$ line running through that point.
The value shifts from $0$ to $1$ at that line.  However, there is no difficulty in the larger space $D_D$.

\subsection{The Service-Time Distribution as a Mixture}\label{cdfmixture}

The general service-time c.d.f. $F$ has at most countably many discontinuity points.
Let $p_d$ $(p_c)$ be the total probability mass at the discontinuity (continuity) points, i.e.,
$p_d \equiv \sum_{x\geq 0} \Delta F(x) \leq 1$ and $p_c = 1-p_d \leq 1$, where $\Delta F(x) \equiv F(x) - F(x-)$.
To focus on the interesting case, suppose that $0 < p_d < 1$.
We order the discontinuity points by the size of their probability mass
in decreasing order (using the natural order in case of ties); i.e., let $\{\bar{x}_1,  \bar{x}_{2}, ... \}$ be such that
 $\Delta F(\bar{x}_i) \geq \Delta F(\bar{x}_{i+1})$.  Define two proper c.d.f.'s $F_c$ and $F_d$
for a continuous random variable $\eta^c$ and a discrete random variable $\eta^d$, respectively, by
$$
 F_c(x) \equiv P(\eta^c \leq x) \equiv \frac{1}{p_c} \Big(F(x) - \sum_{y\leq x} \Delta F(y) \Big), \quad x \geq 0,
$$
and
$$
F_d(x) \equiv \sum_{j: \bar{x}_j \leq x} P(\eta^d = \bar{x}_j) ,
\qandq p_{d,i} \equiv P(\eta^d = \bar{x}_i) \equiv \frac{\Delta F(\bar{x}_i)}{p_d}, \quad x\geq 0.
$$
Note that $F$ can be represented as the mixture $F = p_c F_c + p_d F_d$.

Let $A^c_{n}(t)$, $A^d_{n}(t)$ and  $A^d_{n, i}(t)$
count the number of arrivals by time $t$ with continuous service time,
with a discrete service time, and with a deterministic service time $\bar{x}_i$, $i =1,2,...$, respectively.
Clearly, $A^d_{n}(t) = \sum_{i=1}^{\infty} A^d_{n,i} (t)$ and $A_n (t) = A^d_{n}(t) + A^c_{n}(t)$ for $t\geq 0$. Define
the LLN-scaled processes $ \bar{A}_n^c \equiv n^{-1}A_n^c$, $ \bar{A}_n^d \equiv n^{-1}A_n^d$, and $ \bar{A}_{n,i}^d \equiv n^{-1} A_{n,i}^d$.

Under Assumptions 1 and 2, for a general service-time c.d.f.,
we can decompose the system into two subsystems,
one with arrival processes $A_n^c$ and service-time distribution $F_c$
and the other with arrival processes $A^d_n $ and discrete service times $\{\bar{x}_i: i\geq 1\}$ with distribution $F_d$.
We will adopt the method in \cite{KP97} to analyze the first subsystem in the space $D_D$,
then the method in \cite{GW91} to analyze the second subsystem,
and then we put them together to obtain the limits for the whole system.

\section{Main Results}\label{secMain}

In this section, we state the main results of this paper:
 the FWLLN and FCLT for the scaled processes associated with $Q^r_n$ and $W^r_n$,
 along with the closely related processes.  We give the proofs in \S \ref{secProof}.
Define the LLN-scaled processes $\bar{Q}^r_n \equiv \{\bar{Q}^r_n(t,y), t\geq 0, y\geq 0\}$
by
\beql{c1}
\bar{Q}^r_n(t,y) \equiv \frac{Q^r_n (t,y)}{n},
\eeq
and similarly for the processes $\bar{Q}^e_n$, $\bar{Q}^t_n$, $\bar{D}_n$, $\bar{W}^r_n$, $\bar{W}^t_n$, $\bar{I}_n$ and $\bar{C}_n$.
Define the LLN-scaled processes $\bar{F}^e_n \equiv \{\bar{F}^e_n(t,y), t\geq 0, 0 \leq y \leq  t\}$ and $\bar{F}^{r,c}_n  \equiv \{\bar{F}^{r,c}_n(t,y), t\geq 0, y\geq 0\}$ by
\beql{Fage}
\bar{F}^e_n (t,y) \equiv \bar{Q}^e_n (t,y)/\bar{Q}^t_n (t) \qandq  \bar{F}^{r,c}_n (t,y) \equiv \bar{Q}^r_n (t,y)/\bar{Q}^t_n (t),
\eeq
where $\bar{F}^e_n (t,y)$ and $\bar{F}^{r,c}_n (t,y)$ are defined to be $0$ if $\bar{Q}^t_n (t) = 0$ for some $t$.
By Lemma \ref{repQn},
\beql{E:bQnM}
\bar{Q}^r_n(t,y) =
  \int_0^t F^c (t+y-s) d\bar{a}(s) + \frac{1}{ \sqrt{n}}(\hat{X}_{n,1}(t,y) + \hat{X}_{n,2}(t,y)),\quad t,y \geq 0.
\eeq

When we focus on the amount of work, as in the workload processes,
we use the {\em stationary-excess} (or residual-lifetime) c.d.f.
associated with the service-time c.d.f. $F$ (assumed to have finite mean $\mu^{-1}$), defined by
\beql{se}
F_e (x) \equiv \mu \int_{0}^{x} F^c (s) \, d s, \quad x \ge 0.
\eeq
The mean of $F_e$ is $ E[\eta^2]/2E[\eta] = (c^2_s + 1)/2\mu$; that will be used in part (c) of Theorem \ref{QFWLLN} below.

\begin{theorem} {\em $($FWLLN$)$}\label{QFWLLN}

$($a$)$ Under Assumptions 1 and 2,
\begin{eqnarray}\label{E:fwlln}
&&\left(\bar{A}_n, \bar{A}_n^c, \bar{A}_n^d, \{ \bar{A}_{n,i}^d: i\geq 1\}, \bar{K}_n,\bar{Q}^r_n, \bar{Q}^t_n,\bar{Q}^e_n, \bar{F}^e_n, \bar{F}^{r,c}_n, \bar{D}_n \right) \nonumber \\
&& \Rightarrow \left(\bar{a},\bar{a}^c, \bar{a}^d, \{\bar{a}^d_i: i\geq 1\},  \bar{k},\bar{q}^r,\bar{q}^t,\bar{q}^e,\bar{f}^e, \bar{f}^{r,c}, \bar{d} \right)
\end{eqnarray}
in $D^3 \times D^{\infty} \times D_D^2 \times D \times D_D^3 \times D$ as $n \ra \infty$ w.p.1,
  where the limits are deterministic functions:  $\bar{a}$ is the limit in {\em \eqn{b1}},  $\bar{a}^c \equiv p_c \bar{a}$, $\bar{a}^d \equiv p_d \bar{a}$, $\bar{a}^d_i \equiv p_{d,i} \bar{a}^d$, for $i \geq 1$, $\bar{k}(t,x) \equiv tF(x)$ in {\em \eqn{b7}},
  \begin{eqnarray}\label{E:q}
\bar{q}^r(t,y) &\equiv&  \int_0^t F^c (t+y-s) d \bar{a}(s), \quad t \ge 0, \quad y \ge 0,
\end{eqnarray}
  \begin{eqnarray}\label{E:qe}
\bar{q}^e(t,y) &\equiv&  \int_{t-y}^t F^c (t-s) d \bar{a}(s), \quad t \ge 0, \quad 0 \leq y \leq t,
\end{eqnarray}
$\bar{q}^t (t) \equiv \bar{q}^r (t, 0) = \bar{q}^e (t, t)$,
$\bar{f}^e (t,y) \equiv \bar{q}^e (t,y)/\bar{q}^t (t)$,
$\bar{f}^{r,c} (t,y) \equiv \bar{q}^{r} (t,y)/\bar{q}^t (t)$ and $\bar{d} = \bar{a} - \bar{q}^t$.

$($b$)$  If, in addition to the assumptions in part $($a$)$, $\bar{a} (t) = \lambda t$, $t \ge 0$,
 and the service-time c.d.f. $F$
has finite mean $\mu^{-1}$, then
\beql{fwlln2}
\left(\bar{W}^r_n, \bar{W}^t_n, \bar{I}_n, \bar{C}_n \right) \Rightarrow \left(\bar{w}^r,\bar{w}^t, \bar{i}, \bar{c} \right)
\qinq D_D \times D^3 \qasq n \ra \infty \quad w.p.1,
\eeq
jointly with the limits in {\em \eqn{E:fwlln}},
where
\begin{eqnarray}\label{E:v}
\bar{w}^r(t,y) & \equiv & \lambda \int_y^{\infty} \bar{q}^r (t,x) dx, \quad t \ge 0, \quad y \ge 0, \nonumber \\
&  &  = \lambda \int_y^{\infty} \Big( \int_0^t F^c (t+x-s) d s \Big) \, dx
  = \frac{\lambda}{ \mu} \int_0^t F_e ^c (y+s) d s, \nonumber \\
\bar{w}^t (t) & \equiv & \bar{w}^r (t,0) = \frac{\lambda}{ \mu} \int_0^t F_e ^c (s) d s,  \nonumber \\
\bar{i} (t) & \equiv & \frac{\lambda t}{\mu} \qandq
 \bar{c} (t)  \equiv  \bar{i} (t)  - \bar{w^t} (t)   =\frac{\lambda}{ \mu} \int_0^t F_e (s) d s ,
\end{eqnarray}
for $F_e$ in {\em \eqn{se}}.

$($c$)$ If, in addition to the assumptions of parts $($a$)$ and $($b$)$, $E[\eta^2] < \infty$, then
\beql{b23}
\bar{w}^r(t,y)  \ra \frac{\lambda}{ \mu} \int_0^\infty F_e ^c (y+s) d s < \infty
\qandq \bar{w}^t (t) \ra  \frac{\lambda (c^2_s + 1)}{2 \mu^2} \qasq t\ra\infty.
\eeq
\end{theorem}

We obtain Theorem \ref{QFWLLN} as an immediate corollary to the following FCLT,
which exploits centering by the deterministic limits above.
For the FCLT, define the normalized processes
\beql{b24}
\hat{Q}^r_n (t,y)  \equiv  \sqrt{n} (\bar{Q}^r_n(t,y)  - \bar{q}^r (t,y)),
\eeq
for $t \ge 0$ and $y \ge 0$, and similarly for the other processes, using the centering terms above.
By \eqn{E:bQnM} and \eqn{E:q},
\beql{E:hQnM}
\hat{Q}^r_n(t,y) = \hat{X}^r_{n,1}(t,y) + \hat{X}^r_{n,2}(t,y), \quad t \ge 0, \quad y \ge 0.
\eeq
Moreover,
\begin{eqnarray*}
\hat{F}_n^{r,c}(t,y) &\equiv& \sqrt{n} (\bar{F}^{r,c}_n(t,y)  - \bar{f}^{r,c} (t,y)) \\
&=& \bar{Q}_n^t(t)^{-1} \big( \hat{Q}^r_n(t,y)  - \hat{Q}^t_n(t) \bar{f}^{r,c}(t,y)  \big), \quad t\geq 0, \quad y \geq 0,
\end{eqnarray*}
and
\begin{eqnarray*}
\hat{F}_n^e(t,y) &\equiv& \sqrt{n} (\bar{F}^e_n(t,y)  - \bar{f}^e (t,y)) \\
&=& \bar{Q}_n^t(t)^{-1} \big( \hat{Q}^e_n(t,y)  - \hat{Q}^t_n(t) \bar{f}^e(t,y)  \big), \quad t\geq 0, \quad 0 \leq y \leq t.
\end{eqnarray*}

Define the CLT-scaled processes $\hat{A}_n^c \equiv \{\hat{A}_n^c(t): t\geq 0\}$,
$\hat{A}_n^d \equiv \{\hat{A}_n^d(t): t\geq 0\}$ and $\hat{A}_{n,i}^d \equiv \{\hat{A}_{n,i}^d(t): t\geq 0\}$ by
\begin{eqnarray*}
  \hat{A}_n^c(t) \equiv n^{1/2}(\bar{A}_n^c(t) - \bar{a}^c(t)),
\quad \hat{A}_n^d(t) \equiv n^{1/2}( \bar{A}_n^d(t) - \bar{a}^d(t)), \quad \hat{A}_{n,i}^d(t) \equiv n^{1/2}(\bar{A}_{n,i}^d(t) - \bar{a}^d_i(t)),
\end{eqnarray*}
for  $t \geq 0$ and $i \geq 1$.

The joint deterministic limits in Theorem \ref{QFWLLN} are equivalent to the separate one-dimensional limits,
but that is not true for the FCLT generalization below.  Let $\circ$ be the composition function, i.e., $(x \circ y) (t) \equiv x(y (t))$,
$t \ge 0$.  Let $\deq$ mean equality in distribution.

\begin{theorem} {\em $($FCLT$)$}\label{QFCLT}

$($a$)$ Under Assumptions 1 and 2,
\begin{eqnarray} \label{E:fclt}
&&(\hat{A}_n, \hat{A}_n^c, \hat{A}_n^d, \{ \hat{A}_{n,i}^d: i\geq 1\}, \hat{K}_n, \hat{Q}^r_n, \hat{Q}^t_n, \hat{Q}^e_n, \hat{F}_n^{r,c}, \hat{F}_n^e, \hat{D}_n) \nonumber \\
&& \Rightarrow
(\hat{A}, \hat{A}^c, \hat{A}^d, \{\hat{A}^d_i: i\geq 1\},  \hat{K}, \hat{Q}^r, \hat{Q}^t, \hat{Q}^e, \hat{F}^{r,c}, \hat{F}^e, \hat{D})
\end{eqnarray}
in $D^3 \times D^{\infty} \times D_D^2 \times D \times D_D^3 \times D$ as $n \ra \infty$,
 where $\hat{A}$ is the limit in {\em \eqn{b2}},
  \begin{eqnarray}\label{E:Asplit}
\hat{A}^c & = & p_c \hat{A} + S^c \circ \bar{a}, \quad \hat{A}^d = p_d \hat{A} +S^d \circ \bar{a} ,
\quad \hat{A}^d_i = p_d p_{d,i} \hat{A} +S^{d}_{i} \circ \bar{a}, \nonumber \\
S^c &= & - S^d, \quad S^c \deq \sqrt{p_c(1-p_c)} B, \quad S^d \deq \sqrt{p_d(1-p_d)} B, \nonumber \\
S^d_i &\deq & \sqrt{p_d p_{d,i}(1-p_d p_{d,i})} B, \quad i \geq 1,
\end{eqnarray}
where $B$ is a standard BM, independent of $\hat{A}$, and
the process $(S^c, S^d, \{S^{d}_i: i \geq 1\})$
is an infinite-dimensional BM with mean $0$ and covariance matrix $\bC$ where $\bC_{c,c} = p_c(1-p_c)$,
$\bC_{d,d} = p_d(1-p_d)$, $\bC_{c,d} = \bC_{d,c}=  - p_c p_d$, $\bC_{i,i} = p_dp_{d,i}(1-p_d p_{d,i})$
for $i \geq 1$, $\bC_{i,c}= \bC_{c,i} = - p_c p_dp_{d,i} $, $\bC_{i,d}= \bC_{d,i} = - p_d^2 p_{d,i} $
and $\bC_{i,j} =  - p_d^2 p_{d,i} p_{d,j} $ for $i \neq j$,
 and the representations for $\hat{Q}^r$ and $\hat{Q}^e$ are
 \begin{eqnarray}\label{E:hQrg}
\hat{Q}^r(t,y) &=& \hat{X}^{c,r}_1 (t,y) + \hat{X}^{c,r}_2 (t,y) + \hat{X}^{d,r}(t,y), \quad t\geq 0, \quad y\geq 0, \\
\hat{Q}^e(t,y) &=& \hat{X}^{c,e}_{1} (t,y) + \hat{X}^{c,e}_{2} (t,y) + \hat{X}^{d,e}(t,y), \quad t\geq 0, \quad 0 \leq y \leq t, \nonumber
\end{eqnarray}
where
\begin{eqnarray} \label{E:X1c}
\hat{X}^{c,r}_1(t,y) &\equiv &  \int_0^t F_c^c (t+y-s) d \hat{A}^c(s), \quad \hat{X}^{c,e}_1(t,y) \equiv   \int_{t-y}^t F_c^c (t-s) d \hat{A}^c(s), \\
\hat{X}^{c,r}_2(t,y) & \equiv &  \int_0^t\int_0^{\infty}\mathbf{1}(s+x > t+y)\, d\hat{K}^c(\bar{a}^c(s),x),   \nonumber\\
\hat{X}^{c,e}_2(t,y) & \equiv &  \int_{t-y}^t\int_0^{t}\mathbf{1}(s+x > t)\, d\hat{K}^c(\bar{a}^c(s),x),  \nonumber \\
\hat{X}^{d,r}(t,y) &\equiv& \sum_{i=1}^{\infty} (\hat{A}^d_i(t) - \hat{A}^d_i(t - (\bar{x}_i - y)^{+}) ),  \nonumber \\
\hat{X}^{d,e}(t,y) &\equiv& \sum_{i=1}^{\infty} (\hat{A}^d_i(t) - \hat{A}^d_i(t - (\bar{x}_i  \wedge y)) ),  \nonumber
\end{eqnarray}
with $\hat{K}^c(\bar{a}^c(s),x) = U(\bar{a}^c(s), F_c(x))$, which is independent of $\hat{A}$.
 $\hat{Q}^t (t) \equiv \hat{Q}^r (t,0)$,
 $\hat{Q}^e (t,y) \equiv \hat{Q}^t (t) - \hat{Q}^r (t-y,y)$,
 $\hat{F}^{r,c}(t,y) \equiv \bar{q}^t(t)^{-1} (\hat{Q}^{r}(t,y) - \hat{Q}^t(t) f^{r,c}(t,y) )$,
  $\hat{F}^e(t,y) \equiv \bar{q}^t(t)^{-1} (\hat{Q}^e(t,y) - \hat{Q}^t(t) f^e(t,y) )$,
and $\hat{D} = \hat{A} - \hat{Q}^t$.
 All these limit processes are continuous.   If, in addition,
$\hat{A} =  B_a \circ \bar{a}$, as when $A_n$ is nonhomogeneous Poisson,
then $\hat{A}^d$ and $\hat{A}^c$ are independent, and thus $\hat{X}^{c,r}_1$, $\hat{X}^{c,r}_2$ and $\hat{X}^{d,r}$ are mutually independent.

$($b$)$ If, in addition to the assumptions in part $($a$)$, $\bar{a} (t) = \lambda t$, $t \ge 0$,
 and the service-time c.d.f. $F$
has finite mean $\mu^{-1}$, then
$(\hat{W}^r_n, \hat{W}^t_n) \Rightarrow (\hat{W}^r,\hat{W}^t)$ in $D_D \times D$ as $n \ra \infty$
jointly with the limits in {\em \eqn{E:fclt}},
where
\beql{E:W}
\hat{W}^r (t,y)  \equiv  \int_y^{\infty}  \hat{Q}^r (t,x) \, dx,
\qandq \hat{W}^t (t)   \equiv  \hat{W}^r (t,0) = \int_0^{\infty}  \hat{Q}^r (t,x) \, dx.
\eeq

$($c$)$ If, in addition to the assumptions in parts $($a$)$ and $($b$)$, $E[\eta^2] < \infty$, then
$(\hat{I}_n, \hat{C}_n) \Rightarrow (\hat{I}, \hat{C})$ in $D^2$ as $n \ra \infty$
jointly with the limits above,
where
\beql{E:IC}
\hat{I} (t)     \equiv  \sqrt{\lambda c^2_s} B_s (t) + \mu^{-1} \hat{A} \qandq
\hat{C} (t)     \equiv  \hat{I} (t) - \hat{W}^t (t), \quad t \ge 0,
\eeq
with $B_s$ being a standard BM independent of $\hat{A}$.
\end{theorem}

\begin{remark}\label{rm1}{\em
   The limit processes $\hat{Q}^r$ and  $\hat{Q}^e$ can also be expressed as the sum of the
following three mutually independent processes
\begin{eqnarray}\label{E:hQrgs}
\hat{Q}^r(t,y) &=& \hat{X}^r_1 (t,y) + \hat{X}^{c,r}_2 (t,y) + \hat{X}^r_{3}(t,y), \quad t\geq 0, \quad y\geq 0, \\
\hat{Q}^e(t,y) &=& \hat{X}^e_1 (t,y) + \hat{X}^{c,e}_2 (t,y) + \hat{X}^e_3(t,y), \quad t\geq 0, \quad 0 \leq y \leq t, \nonumber
\end{eqnarray}
where
\begin{equation} \label{E:Xr1}
\hat{X}^r_1(t,y) \equiv \int_0^t F^c (t+y-s) d \hat{A}(s), \quad
\hat{X}^e_1(t,y)  \equiv  \int_{t-y}^t F^c (t-s) d \hat{A}(s),
\end{equation}
 \begin{eqnarray*}
 \hat{X}^r_3(t,y) &\equiv& \int_0^t F_c^c(t+y-s) d S^c (\bar{a}(s)) + \sum_{i=1}^{\infty}  \big( S^d_i( \bar{a} (t) )
 - S^d_i ( \bar{a}(t - (\bar{x}_i - y)^{+}) ) \big), \\
  \hat{X}^e_3(t,y) &\equiv& \int_{t-y}^t F_c^c(t-s) d S^c (\bar{a}(s)) + \sum_{i=1}^{\infty}  \big( S^d_i( \bar{a} (t) )
 - S^d_i ( \bar{a}(t - (\bar{x}_i \wedge y)) ) \big).
  \end{eqnarray*}
The asymptotic variability of the arrival process is captured by
$\hat{A}$, which appears only in $\hat{X}^r_{1}$ and $\hat{X}^e_{1}$; the asymptotic variability of the service process is captured
by $\hat{K}^c$, which appears only in $\hat{X}^{c,r}_{2}$ and $\hat{X}^{c,e}_{2}$; while the asymptotic variability of service-time
splitting is captured by $S^c$ and $S^d_i$, which appears only in $\hat{X}^r_{3}$ and $\hat{X}^e_{3}$.  Thus, in some sense,
there is additivity of stochastic effects,
as pointed out in \cite{L88, KP97}, but this might be misinterpreted.  Notice that {\em both} $\hat{X}^r_{1}$
and $\hat{X}^r_{2}$ depend on the full service-time c.d.f.  $F$, not just its mean.  On the other hand, the
arrival process beyond its deterministic rate only appears in $\hat{X}^r_{1}$ and $\hat{X}^e_{1}$,
so that there is a genuine asymptotic insensitivity to the arrival
process beyond its rate in $\hat{X}^{c,r}_{2}$ and $\hat{X}^{c,e}_{2}$.

This claim holds because, by \eqn{E:Asplit}, we can write $ \hat{X}^c_1(t,y)$ and  $\hat{X}^d(t,y) $ in \eqn{E:hQrg} as
\begin{eqnarray*}
\hat{X}^{c,r}_1(t,y)  & = &  \int_0^t F_c^c (t+y-s) d ( p_c \hat{A}(t) + S^c (\bar{a}(s)) ) \\
&=& \int_0^t \Big(F^c(t+y-s) -\sum_{u > t+y-s} \Delta F(u)\Big) d  (\hat{A}(t) + p_c^{-1}  S^c (\bar{a}(s)) ),
\end{eqnarray*}
and
\begin{eqnarray*}
  \hat{X}^{d,r}(t,y)
 & = & \sum_{i=1}^{\infty} \Big[   p_d p_{d,i} \big( \hat{A} (t)
 - \hat{A} (t - (\bar{x}_i - y)^{+}) ) \big)  + \big( S^d_i( \bar{a}( t) )
 - S^d_i (\bar{a} (t - (\bar{x}_i - y)^{+})) \big) \Big] \\
 &=&  \int_0^t \Big(\sum_{u > t+y-s} \Delta F(u)\Big) d  \hat{A}(t)
 +  \sum_{i=1}^{\infty}  \big( S^d_i( \bar{a}( t) ) - S^d_i (\bar{a} (t - (\bar{x}_i - y)^{+})) \big),
\end{eqnarray*}
which implies that $\hat{X}^{c,r}_1(t,y) + \hat{X}^{d,r}(t,y)  = \hat{X}^r_1(t,y) + \hat{X}^r_3(t,y)$ for each $t\geq 0$ and $y\geq 0$.
Similarly,  $\hat{X}^{c,e}_1(t,y) + \hat{X}^{d,e}(t,y)  = \hat{X}^e_1(t,y) + \hat{X}^e_3(t,y)$ holds. ~~~\bsq
}
\end{remark}

\begin{remark}\label{rm2}
{\em The two integrals in the expression for $\hat{Q}^r$
are stochastic integrals.  The first integral for $\hat{X}^{c,r}_1$ (or $\hat{X}^r_1$) is a standard Ito integral if $\hat{A}$ is a (time-changed) Brownian motion;
otherwise, the expression for  $\hat{X}^{c,r}_1$  (or $\hat{X}^r_1$) is interpreted as the form after integration by parts.
The relevant version of integration by parts for $\hat{X}^r_{n,1}$ and $\hat{X}^r_{1}$
is given in Bremaud \cite{B81}, p.336.  For $\hat{X}^r_{n,1}$, it yields
\beql{Xn1parts}
\hat{X}^r_{n,1} (t,y)  =   F^c (y) \hat{A}_n(t) - \int_0^t  \hat{A}_n(s-) \,d F(t + y -s),
\eeq
and similarly for $\hat{X}_{1}$.  The left limit $\hat{A}_n(s-)$ in \eqn{Xn1parts} is only needed if the
functions $F$ and $\hat{A}_n$ have common discontinuities with positive probability.
The second integral for $\hat{X}_2$ is either understood as the stochastic integrals with respect to
two-parameter processes of the first type, or in the mean-square sense, as in \cite{KP97}; see \S \ref{secFiDi}.

In the literature, several types of stochastic integrals with respect to two-parameter processes have been defined.
The first type of integral was first defined for
two-parameter Brownian sheets by Cairoli \cite{Cair} (see also \cite{Walsh}),  generalizing the definition of Ito's integral directly. It was
generalized to $n$-parameter Brownian sheets by Wong and Zakai \cite{WZ74} and to general martingales by Cairoli and Walsh \cite{CW75}.
Even more generalization appears in Wong and Zakai \cite{WZ77}. We refer to Koshnevisan \cite{K02} for a relatively complete review. The important property we apply here is the isometry property,
analogous to the Ito isometry property.
~~~\bsq
}
\end{remark}

\begin{remark}\label{rm3}
{\em
We remark that if the service-time c.d.f. $F$ is discontinuous, the process $\hat{X}_2$ defined by
$$
\hat{X}_2(t,y)  \equiv   \int_0^t\int_0^{\infty}\mathbf{1}(s+x > t+y)\, d\hat{K}(\bar{a}(s),x)
$$
 is only continuous in $t$,
but not in $y$, and in fact, it is not even in the space $D_D$.
The continuity of $\hat{X}_2$ and $\hat{Q}^t$ in $t$ can be obtained as in Lemma 5.1 of \cite{KP97}.
To see that $\hat{X}_2$ need not be in $D_D$,
suppose that $F$ is the mixture of two point masses $y_1>0$ and $y_2>0$.
Then, applying \eqn{Kinc} below, we see that, for each $t\geq 0$, $\hat{X}_2(t, y) = 0$ for all $y\geq 0$ except $y_1$ and $y_2$,
so that $\hat{X}_2(t, \cdot) \notin D$.
That property follows from \eqn{Kinc} because $\Delta_{\hat{K}} (t_1, t_2, x_1, x_2) = 0$ for $0 < x_1 < x_2$ unless either
$y_1 < x_1 < y_2$ or $x_1 < y_1 < x_2 < y_2$.  That means that the random measure attaches all mass on the strips $x = y_1$
and $x = y_2$.
Incidentally,
in this example, $\hat{X}_2(t, \cdot)$
is an element of the space $E$ in Chapter 15 of \cite{W02}.
That explains why we split the general distribution into a mixture of a discrete distribution and a continuous distribution.
~~~\bsq
}
\end{remark}

We now establish additional results in the standard case for the fluid limits.
In particular, we will obtain an analog of the classic result
for the $M/GI/\infty$ model, stating that in steady state both the elapsed service times and the residual service times
are distributed as mutually independent random variables, each with c.d.f. $F_e$ in \eqn{se}.
We will see that the limiting empirical age distribution
is precisely $F_e$, just as is true for the prelimit processes with a Poisson arrival process.

\begin{coro}{\em $($the standard case$)$}\label{corStd}
Consider the standard case in which $\bar{a} (t) = \lambda t$, $t \ge 0$, and $\hat{A} = \sqrt{\lambda c^2_a}B_a$,
where $B_a$ is a standard BM.  Assume that the service-time distribution $F$ has finite mean $\mu^{-1}$.
Under Assumptions 1 and 2, the limits in {\em \eqn{E:fwlln}} hold with
\begin{eqnarray}\label{Std}
\bar{q}^r(t,y) & \equiv & \lambda \int_0^t F^c (t+y-s)\, d s =
\lambda \int_0^t F^c (y + s)\, d s \nonumber {} \\
&&{}{} \ra \left(\lambda/\mu\right) F_e^c (y) \qasq t \ra \infty, \nonumber \\
\bar{q}^e(t,y) & \equiv & \lambda \int_{t-y}^t F^c (t-s)\, d s =
\lambda \int_0^y F^c ( s)\, d s = \left(\lambda/\mu\right) F_e (y),  \qforq t \geq 0, \nonumber \\
\bar{f}^e (t,y) & \equiv  & \bar{q}^e (t,y)/\bar{q}^t (t)  \ra  F_e (y) \qasq t\ra\infty,  \nonumber \\
\bar{f}^{r,c} (t,y) & \equiv  &  \bar{q}^r (t,y)/\bar{q}^t (t)  \ra  F_e^c (y) \qasq t \ra \infty.
\end{eqnarray}
\end{coro}



%



\section{Characterizing the Limit Processes}\label{secChar}

We now show that the two-parameter queue-length limit processes, $\hat{Q}^r (t, y)$ and $\hat{Q}^e (t, y)$,
constitute continuous Brownian analogs of the
Poisson-random-measure representation for the $M/GI/\infty$ model \cite{EMW93}.  (But the limit is only identical to the limit for the
$M/GI/\infty$ model when $c_a^2 = 1$.)
A key role here is played by the \textit{transformed Kiefer process} $\hat{K} (t, x) \equiv U(t, F(x)) = W(t, F(x)) - F(x) W(t,1)$.
Any finite number of $\hat{K}$-increments,
\begin{eqnarray} \label{Kinc}
\Delta_{\hat{K}} (t_1, t_2, x_1, x_2)  &\equiv& \hat{K} (t_2, x_2) - \hat{K} (t_2, x_1) - \hat{K} (t_1, x_2) +  \hat{K} (t_1, x_1) \\
&=& \Delta_W(t_1,t_2,F(x_1),F(x_2)) - (F(x_2) - F(x_1))(W(t_2,1) - W(t_1,1)) \nonumber
\end{eqnarray}
for $0 \le t_1 < t_2$ and $0 \le x_1 < x_2$,
are independent random variables provided that the rectangles $(t_1,t_2]\times (x_1,x_2]$ have disjoint horizontal time intervals $(t_1, t_2]$.

\begin{figure}[h!]
\centering
   \includegraphics[width=4.5in,height=2.8in]{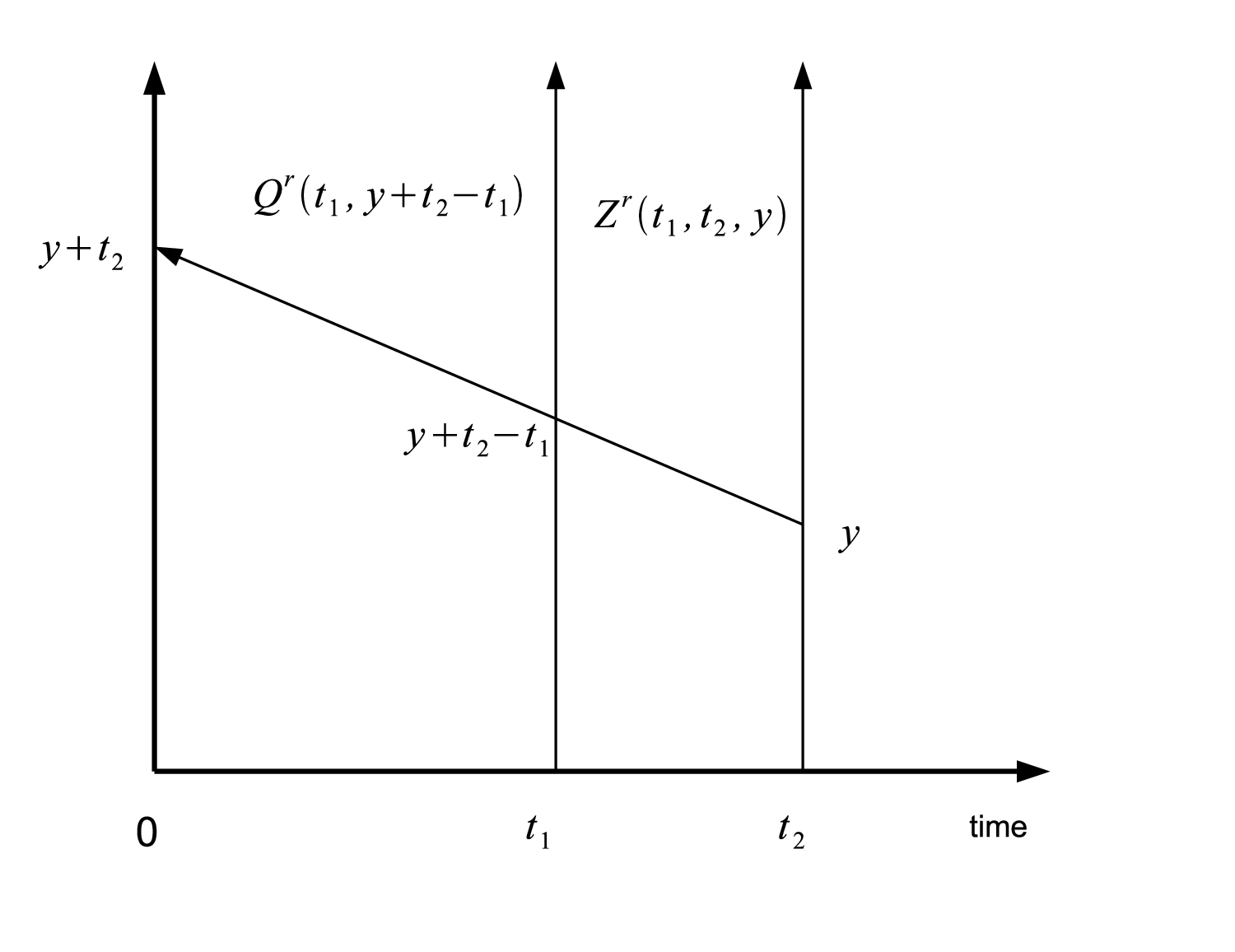}
\caption{The basic decomposition for $Q^r (t,y)$.}
\label{fig2}
\end{figure}

We only treat $\hat{Q}^r$ here.
If the limit process $\hat{A}$ has independent increments, then so does $\hat{Q}^r$,
 provided that it is viewed as a function-valued process with the argument $t$.
The limit processes $\hat{Q}^r$ is then a Markov process in $D_D$ (only considering the argument $t$).
This result can be based on a
basic decomposition, depicted in Figure \ref{fig2}.

\begin{theorem} {\em $($decompositions, independent increments and the Markov property for $\hat{Q}^r)$}\label{thMarkov}
The limiting random variables $\hat{X}^{c,r}_1 (t,y)$, $\hat{X}^{c,r}_2 (t,y)$, $\hat{X}^{d,r} (t,y)$ and $\hat{Q}^r (t,y)$ in Theorem {\em \ref{QFCLT}} admit the decompositions
\begin{eqnarray}\label{E:decomp}
\hat{X}^{c,r}_i (t_2,y) & = & \hat{X}^{c,r}_i (t_1,y+t_2-t_1) + Z^{c,r}_i (t_1,t_2,y), \qforq i = 1,2, \qandq  t_2 >t_1\geq 0, \nonumber \\
\hat{X}^{d,r} (t_2,y) &=& \hat{X}^{d,r} (t_1,y+t_2-t_1) + Z^{d,r} (t_1,t_2,y),  \quad t_2 >t_1\geq 0, \nonumber\\
\hat{Q}^r (t_2,y) & = & \hat{Q}^r (t_1,y+t_2-t_1) + Z^r (t_1,t_2,y), \quad t_2 >t_1\geq 0,
\end{eqnarray}
where $y \ge 0$, $Z^r \equiv Z^{c,r}_1 + Z^{c,r}_2 + Z^{d,r}$, and
\begin{eqnarray*}\label{E:integrals}
Z^{c,r}_1 (t_1,t_2,y) & \equiv &  \int_{t_1}^{t_2} F_c^c (t+y-s) d \hat{A}^c(s), \nonumber \\
Z^{c,r}_2 (t_1,t_2,y) & \equiv &  \int_{t_1}^{t_2} \int_{0}^{\infty}\mathbf{1}(s+x > t+y)\, d \hat{K}^c(\bar{a}^c(s),x), \nonumber\\
Z^{d,r} (t_1,t_2,y) &\equiv& \sum_{i=1}^{\infty} \big[ (\hat{A}^d_i(t_2) -\hat{A}^d_i(t_1)  ) -  ( \hat{A}^d_i(t_2 - (\bar{x}_i - y)^{+}) - \hat{A}^d_i(t_1 - (\bar{x}_i - y)^{+}) ) \big].
\end{eqnarray*}
If, in addition to the assumptions of Theorem {\em \ref{QFCLT}}, the limit process
$\hat{A}$ has independent increments, which occurs in the standard case of Corollary {\em \ref{corStd}}, where
$\hat{A}$ is a BM,
then the two random variables on the right in {\em \eqn{E:decomp}} are independent in each case.
Moreover, the four processes
$\{\hat{X}^{c,r}_1 (t,\cdot): t \ge 0\}$, $\{\hat{X}^{c,r}_2 (t,\cdot): t \ge 0\}$, $\{\hat{X}^{d,r} (t,\cdot): t \ge 0\}$  and $\{\hat{Q}^r (t,\cdot): t \ge 0\}$ all have independent
increments, and are thus Markov processes $($with respect to the argument $t$$)$.

\end{theorem}

\paragraph{Proof.}  The decomposition for $\hat{X}^{c,r}_1 (t,y)$, $\hat{X}^{c,r}_2 (t,y)$, $\hat{X}^{d,r}(t,y)$ and $\hat{Q}^r (t,y)$ in \eqn{E:decomp}
is by direct construction, as in Figure \ref{fig2}.
The independent-increments property is inherited from $\hat{K}^c$, $\hat{A}^c$ and $\hat{A}^d$. ~~~\bsq

We now show that the limit processes are Gaussian if $\hat{A}$ is Gaussian, which again is the case if
$\hat{A}$ is BM.  For nonstationary non-Poisson arrival processes ($G_t$), we can construct such $G_t$ processes
(or just think of them) by letting the original arrival processes $\{A_n (t): t \ge 0\}$
be defined by
$A_n (t) \equiv \tilde{A} (n \bar{a} (t))$, $t \ge 0$,
 where $\tilde{A} \equiv \{\tilde{A} (t): t \ge 0\}$
is a rate-$1$ stationary (or asymptotically stationary) stochastic point process, such that $\tilde{A}$ satisfies a FCLT
with limit $\sqrt{c_a^2} B_a$, where $B_a$ is a standard BM.  As a consequence, a natural Gaussian limit process is
$\hat{A} (t) \equiv \sqrt{c_a^2} B_a (\bar{a} (t))$, $t \ge 0$.  Indeed, this occurs for the familiar $M_t$ case,
for which $c_a^2 = 1$.

\begin{theorem} {\em $($Gaussian property$)$}\label{thGaussian}
If, in addition to the assumptions of Theorem {\em \ref{QFCLT}}, the limit process
$\hat{A}$ is Gaussian, then
the limit processes
$\hat{Q}^t, \hat{Q}^e, \hat{Q}^r, \hat{D}, \hat{V}^r, \hat{V}^t $
in {\em \eqn{E:fclt}} are all continuous Gaussian processes.
If $\hat{A} (t) = \sqrt{c^2_a} B_a (\bar{a} (t))$ for $t \ge 0$, where $B_a$ is a standard BM,  then
for each fixed $t\geq 0$ and $y\geq 0$,
\begin{eqnarray} \label{E:varhQg}
 \hat{Q}^r(t,y) \deq N (0, \sigma_{q,r}^2(t,y)), \quad \hat{Q}^e(t,y) \deq N (0, \sigma_{q,e}^2(t,y)), \quad  \hat{W}^r(t,y) \deq N (0, \sigma_{w}^2(t,y)),
\end{eqnarray}
where
\begin{eqnarray*}
\sigma_{q,r}^2(t,y) &=& (c_a^2-1) \int_0^t F^c(t+y-s)^2 d \bar{a}(s) +  \int_0^t  F^c(t+y-s) d \bar{a}(s),\\
 \sigma_{q,e}^2(t,y) &=&  (c_a^2-1) \int_{t-y}^t F^c(t-s)^2 d \bar{a}(s) + \int_{t-y}^t   F^c(t-s) d \bar{a}(s), \\
 \sigma_{w}^2(t,y) &=&  c_a^2 \int_y^{\infty} \int_y^{\infty} \int_0^t F^c(t+x-s) F^c(t+z-s) d \bar{a}(s) dx dz {}\\
 && +  \int_y^{\infty} \int_y^{\infty} \int_0^t  F(t + x\wedge z -s) F^c(t+x\vee z -s) d \bar{a}(s) dx dz.
\end{eqnarray*}

\end{theorem}

\paragraph{Proof.}
It is obvious that the limit processes are Gaussian when the limit arrival process $\hat{A}$ is Gaussian. We only need to derive the variance formulas.  We will use \eqn{E:hQrgs} to calculate them and  the mutual independence between the three terms in the expression of $\hat{Q}^r$
gives $\sigma_{q,r}^2(t,y)  = \sigma_{1,r}^2(t,y) + \sigma_{2,c,r}^2(t,y) + \sigma_{3,r}^2(t,y)$, where $\sigma_{1,r}^2(t,y)  = E[(\hat{X}^{r}_1 (t,y) )^2]$, $\sigma_{2,c,r}^2(t,y)  = E[(\hat{X}^{c,r}_2 (t,y) )^2]$ and $\sigma_{3,r}^2(t,y) = E[( \hat{X}^{r}_3(t,y))^2]$.
 By Ito's isometry, we have
$$
\sigma_{1,r}^2(t,y)  = c_a^2 \int_0^t F^c(t+y-s)^2 d \bar{a}(s),
$$
and
\begin{eqnarray*}
 \sigma_{3,r}^2(t,y)
 &=& p_d p_c  \int_0^t  F_c^c(t+y-s)^2 d \bar{a}(s) + \sum_{i=1}^{\infty} \Big( p_dp_{d,i} (1 - p_d p_{d,i}) (\bar{a}(t) - \bar{a}(t - (\bar{x}_i - y)^{+})) \Big){}\\
 && {}- 2  p_d^2 \sum_{i<j}   p_{d,i} p_{d,j}\big( \bar{a}(t) - \bar{a}(t - ((\bar{x}_i \wedge \bar{x}_j) - y)^{+})  \big) {}\\
 &&{} - 2 \sum_{i=1}^{\infty} p_c p_d p_{d,i}\int_0^t F_c^c(t+y-s) d  (\bar{a}(s) - \bar{a}(s- (\bar{x}_i - y)^{+})).
\end{eqnarray*}

Having $\hat{X}^{c,r}_2$ well-defined with continuous paths follows from the definition of stochastic integral
with respect to the Brownian sheet of the first type.  It clearly has mean 0.  Its variance is given by
\begin{eqnarray*}
 \sigma_{2,c,r}^2(t,y)
 & = &
  E\Big[\Big(\int_0^t\int_0^{\infty}\mathbf{1}(s+x> t+y)
d U(\bar{a}^c(s),F_c(x))\Big)^2\Big] \nonumber \\
&=& E\Big[\Big(\int_0^t\int_0^{\infty}\mathbf{1}(s+x> t+y)
d (W(\bar{a}^c(s),F_c(x))-F_c(x)W(\bar{a}^c(s),1))\Big)^2\Big] \nonumber \\
&=& \int_0^t\int_0^{\infty}\mathbf{1}(s+x > t+y) d F_c(x) d \bar{a}^c(s) +\int_0^tF_c^c(t+y-s)^2 d \bar{a}^c(s){} \nonumber \\
& & {} - 2 \int_0^t\int_0^{\infty}\mathbf{1}(s+x > t+y)F_c(t+y-s)dF_c(x)d \bar{a}^c(s)\nonumber \\
&=& \int_0^t F_c(t+y-s) F_c^c(t+y-s)d \bar{a}^c(s),
\end{eqnarray*}
where the second equality uses the identity $U(x,y) = W(x,y) - yW(x,1)$, and the third equality uses the isometry
property of the stochastic integral of the first type with respect to two-parameter
Brownian sheets and also the isometry property of the stochastic Ito's integral.

Notice that
\begin{eqnarray*}
&&p_d p_c\int_0^t F_c^c(t+y-s)^2 d  \bar{a}(s)  + \int_0^t F_c(t+y-s) F_c^c(t+y-s)d \bar{a}^c(s) \\
&=& \int_0^t p_c F_c^c(t+y-s) (1- p_c F_c^c(t+y-s) ) d \bar{a}(s).
\end{eqnarray*}
Moreover, $F^c = p_c F_c^c + p_d F_d^c$ and  $F F^c= (1- p_c F_c^c - p_d F_d^c) (p_c F_c^c + p_d F_d^c) = p_c F_c^c(1- p_c F_c^c) +  p_d F_d^c (1- p_d F_d^c) - 2 p_c F_c^c p_d F_d^c$. Then, simple algebra calculation gives the final expression for $\sigma_{q,r}^2(t,y)$. Similar argument applies to the calculation of $\sigma_{q,e}^2(t,y)$.

For the variance of $\hat{W}^r(t,y)$, by the independence of $\hat{X}^{r}_1(t,y)$, $\hat{X}^{c,r}_2(t,y)$ and $\hat{X}^{r}_3(t,y)$, we have
\begin{eqnarray*}
E[\hat{W}^r(t,y)^2]
=   E\Big[ \Big( \int_y^{\infty} \hat{X}^{r}_1(t,x) dx\Big)^2\Big] + E\Big[ \Big( \int_y^{\infty} \hat{X}^{c,r}_2(t,x) dx\Big)^2\Big]  + E\Big[ \Big( \int_y^{\infty} \hat{X}^{r}_3(t,x) dx\Big)^2\Big].
\end{eqnarray*}
Then by an analogous argument, we obtain the variance of $\hat{W}^r(t,y)$. ~~~\bsq

We remark that, for Theorem \ref{thGaussian}, we could also have used an argument analogous to Lemma 5.1 in \cite{KP97}
by understanding the integral in $\hat{X}^{c,r}_2$ as a mean-square limit (\S \ref{secFiDi}).
However, our approach here by applying properties of stochastic integrals of the first type
with respect to two-parameter Brownian sheets  simplifies the proof.
Paralleling the result in Lemma 5.1 \cite{KP97}, we can easily check that for $ 0 \leq t \leq t' $, $ 0 \leq  y \leq y'$,
\begin{eqnarray*}
&& E[(\hat{X}^{c,r}_2(t,y)- \hat{X}^{c,r}_2(t',y'))^2] \\
&=&  \int_0^{t} (F_c(t'+y'-u)-F_c(t+y-u))(1+F_c(t+y-u)-F_c(t'+y'-u)) d\bar{a}^c(u).
\end{eqnarray*}

\begin{coro}  $($\mbox{the special case $c_a^2 = 1$}$)$ \label{coro1Var}
If, in addition to the assumptions of Theorem $\ref{thGaussian}$, $\bar{a}(t) = \int_0^t \lambda(s)ds$ and $c_a^2 = 1$,  then
\beql{E:PoissonCase1}
Var(\hat{Q}^r(t,y)) =
\int_0^{t} F^c(t+y-u) \, \lambda (s) \, d s,
\eeq
for $t\geq 0$ and $y\geq 0$.
The limit $\hat{A}$ and all the
other limits are the same as if the unscaled arrival processes
$\{A_n (t): t \ge 0\}$ are Poisson processes $($possibly
nonhomogeneous$)$.  $($When $A_n$ is Poisson, the prelimit
variables $Q^r_n (t,y)$ and $Q^e_n (t,y)$ are Poisson random
variables for each $t$ and $y$.$)$
 Moreover, as in the Poisson arrival case, for each
$t \ge 0$ and $y \ge 0$, $\hat{Q} ^r (t,y)$ is distributed the
same as the limit of
\begin{equation} \label{E:PoissonCase2}
\hat{\sQ}^r_n(t,y) \equiv \sqrt{n} \Big( \frac{1}{n} \sum_{i = 1}^{Q_{n}^t (t)} \eta_{i} (t,y) - \bar{q}^r(t,y) \Big),
\end{equation}
where $\{\eta_{i} (t,y): i \ge 1\}$ is a sequence of i.i.d. Bernoulli random variables with
\beql{E:PoissonCase3}
P(\eta_{i} (t,y) = 1) = \bar{f}^{r,c}(t,y),
\eeq
which are independent of the total queue length $\hat{Q}^{t}_{n} (t)$.
\end{coro}

\paragraph{Proof.}

We need to justify \eqn{E:PoissonCase2}.  First, we note
that this is the asymptotic generalization of an exact relation
for Poisson arrivals; e.g., see Theorem 2.1 of \cite{GW08}.  Here
we start by defining
$$
\sQ^r_n(t,y) \equiv \sum_{i =1}^{Q_n^t(t)} \eta_i(t,y),
$$
for each $t\geq 0$ and $y\geq 0$.  (In passing, we remark that $\sQ^r_n(t,y) \deq Q_n^r(t,y)$
in the special case of a nonhomogeneous ($M_t$) arrival process, but not more generally.)
By the FWLLN, the fluid scaled processes  $\bar{\sQ}^r_n(t,y)$ converge to the fluid limit $\bar{q}^r(t,y)$ as $n\ra\infty$:
\begin{eqnarray*}
\bar{\sQ}^r_n(t,y) \Rightarrow  \bar{\sQ}^r(t,y)  \equiv E[\eta_i(t,y)] \bar{q}^t(t)= \bar{f}^{r,c}(t,y)\bar{q}^t(t)
 =  \frac{ \bar{q}^r(t,y)}{\bar{q}^t(t) }\bar{q}^t(t) =  \bar{q}^r(t,y).
\end{eqnarray*}
We can write $ \hat{\sQ}^r_n(t,y)$ in \eqn{E:PoissonCase2} as
\begin{eqnarray*}
\hat{\sQ}^r_n(t,y) &=&  \frac{1}{\sqrt{n}}  \sum_{i =1}^{n \bar{Q}_n^t(t)} \big(  \eta_i(t,y) - \bar{f}^{r,c}(t,y) \big)
 + \bar{f}^{r,c}(t,y) \hat{Q}_n^t(t).
\end{eqnarray*}
By FCLT for random walks with i.i.d. increments of mean $0$ and finite variance (Theorem 8.2, \cite{B68}),
continuity of composition in $D$ (Theorem 13.2.2,  \cite{W02}) and Theorems \ref{QFWLLN} and \ref{QFCLT},
we obtain the weak convergence of $ \hat{\sQ}^r_n(t,y)$:
\begin{eqnarray*}
\hat{\sQ}^r_n(t,y)  \Rightarrow \hat{\sQ}^r(t,y) \qinq D_D \qasq n\ra\infty,
\end{eqnarray*}
where
$$
\hat{\sQ}^r(t,y) \equiv \sigma_3(t,y) B_3(\bar{q}^t(t)) + \bar{f}^{r,c}(t,y)  \hat{Q}^t(t)
$$
where $\sigma_3^2(t,y) \equiv   \bar{f}^{r,c}(t,y) (1- \bar{f}^{r,c}(t,y) )$ and $B_3$ is a standard Brownian motion, independent of $\hat{Q}^t(t)$.
Thus, $\hat{\sQ}^r(t,y) $ is Gaussian with mean $0$ and variance
\begin{eqnarray*}
 Var(\hat{\sQ}^r(t,y) )
&=&  \sigma_3^2(t,y)\bar{q}^t(t) + \bar{f}^{r,c}(t,y)^2 Var( \hat{Q}^t(t)) \\
&=&  \bar{f}^{r,c}(t,y) (1- \bar{f}^{r,c}(t,y) )\bar{q}^t(t) + \bar{f}^{r,c}(t,y)^2  \int_0^{t} F^c(t-u) \, \lambda (s) \, d s \\
&=&  \frac{ \bar{q}^r(t,y)}{\bar{q}^t(t) } \Big(1 - \frac{ \bar{q}^r(t,y)}{\bar{q}^t(t) } \Big) \bar{q}^t(t)
+ \frac{ \bar{q}^r(t,y)^2}{\bar{q}^t(t) ^2}\bar{q}^t(t)  \\
&=& \bar{q}^r(t,y)  =  \int_0^{t} F^c(t+y-u) \, \lambda (s) \, d s = Var(\hat{Q}^r(t,y)).
\end{eqnarray*}
Since $\hat{\sQ}^r(t,y) $ and $\hat{Q}^r(t,y) $ are both Gaussian with the same mean and variance, $\hat{\sQ}^r(t,y) $
and $\hat{Q}^r(t,y) $ are equal in distribution.
When the arrival process is $M_t$,
$Q_n^r (t, y)$ has a Poisson distribution for each $n$, $t$ and $y$, so that the variance equals the mean.
  Since $c_a^2 = 1$, the limit must be the same here as in the $M_t$ case.
 ~~~\bsq

We emphasize that Corollary \ref{coro1Var} is consistent with known results for
the $M_t/GI/\infty$ model. The asymptotic equivalence to the
random sum in \eqn{E:PoissonCase2} and \eqn{E:PoissonCase3} is
the asymptotic analog of the property for the $M_t/GI/\infty$
model that, conditional on the number of customers in the system,
the remaining service times are distributed as i.i.d. random
variables with c.d.f. $\bar{f}^{r,c} (t, \cdot)$; e.g., see Theorem 2.1
of \cite{GW08}. This property does not hold for $c_a^2 \not= 1$.

\begin{coro} $($\mbox{the standard case}$)$ \label{coro2Var}
If $\bar{a} (t) = \lambda t$ and $\hat{A} = \sqrt{\lambda c^2_a} B_a$, then the variances of $\hat{Q}^r(t,y)$ and $\hat{Q}^e(t,y)$ are
\begin{eqnarray*}
\sigma_{q,r}^2(t,y)
&=& \lambda (c_a^2-1) \int_0^t F^c(y+s)^2 d s +   \lambda \int_0^t  F^c(y+s) d s \\
&\ra& \lambda (c_a^2-1) \int_y^{\infty} F^c(s)^2 d s +   \lambda \int_{y}^{\infty}  F^c(s) d s \equiv \sigma_{q,r}^2(y)  \qasq t\ra\infty, \quad y\geq 0,
\end{eqnarray*}
and
\begin{eqnarray*}
\sigma_{q,e}^2(t,y)
&=& \lambda (c_a^2-1) \int_0^y F^c(s)^2 d s +   \lambda \int_0^y  F^c(s) d s\equiv \sigma_{q,e}^2(y), \quad t, y\geq 0.
\end{eqnarray*}
Thus, $\hat{Q}^r(t,y) \Rightarrow N(0,\sigma_{q,r}^2(y))$ and $\hat{Q}^e(t,y) \Rightarrow N(0,\sigma_{q,e}^2(y))$ as $t\ra\infty$.
If, in addition, $c_a^2 = 1$, then
\begin{eqnarray*} \label{E:PoissonCase4}
\sigma_{q,r}^2(t,y) &=& \lambda \int_0^t F^c(y+s) d s \ra \lambda \int_y^{\infty} F^c(s) d s  =  \frac{\lambda}{\mu} F_e^c(y), \qasq t\ra\infty, \quad y\geq 0, \\
\sigma_{q,e}^2(t,y) &=& \lambda  \int_0^y F^c(s) d s = \frac{\lambda}{\mu} F_e(y), \quad t, y\geq 0,
\end{eqnarray*}
and $Var(\hat{Q}^t(t)) = \lambda \int_0^t F^c(s) d s \ra \lambda/ \mu$ as $t \ra\infty$.

\end{coro}

\section{Initial Conditions}\label{secInitial}

So far, we considered only new arrivals.  Now we consider customers in the system initially.
Like the generality of the service-time c.d.f., the initial conditions present technical difficulties.
Our assumptions will be similar to those made in \cite{KP97} and
to those for the new arrivals in \S \ref{secPrelim}. However, these assumptions are less realistic here.
Thus, for applications, it is good that the relevance of the initial conditions decreases as time evolves,
because we can think of the system starting in the distant past with just new arrivals, so that we will be able to approximate the
two-parameter processes by the Markov limit processes.

We assume that the remaining service times of
the customers initially in the system are i.i.d., distributed according to some new c.d.f., independent of
the number of customers in the system and everything associated with new arrivals.  That rather strong assumption will actually be justified
if we assume that the initial state we see is the result of an $M_t/GI/\infty$ system, possible with different model parameters,
 that started empty at some previous time.
As noted in Corollary \ref{coro1Var} and the remark before Corollary \ref{coro2Var}, this strong independence property actually holds in
an $M_t/GI/\infty$ model.  Moreover, that representation is
asymptotically correct more generally if $c_a^2 = 1$.  Unfortunately, however, that representation is not asymptotically correct
if $c_a^2 \not= 1$.
Nevertheless, it is a natural candidate approximate initial condition.

Here is our specific framework:
Let $Q^{i,r}_n (y)$ be the number of customers initially in the $n^{\rm th}$ system at time $0$,
not counting new arrivals, who have residual service times strictly greater than $y$.
Let $Q^{i,t}_n \equiv Q^{i,r}_n (0)$ be the total number of customers initially in the $n^{\rm th}$ system and let
$Q^{i,e}_n (y) $ be the number of customers initially in the $n^{\rm th}$ system
that have elapsed service times less than or equal to $y$.
Let $W^{i,r}_n (y)$ and $W^{i,t}_n$ be the corresponding workload processes, defined as in \eqn{E:Vn}.

Let $\bar{Q}^{i,r}_n (y)$ and $\hat{Q}^{i,r}_n (y)$ be the associated scaled processes, defined by
\beql{f1}
\bar{Q}^{i,r}_n (y) \equiv \frac{Q^{i,r}_n (y)}{n} \qandq \hat{Q}^{i,r}_n (y)
\equiv \sqrt{n}(\bar{Q}^{i,r}_n (y) - \bar{q}^{i,r} (y)), \quad y\geq 0,
\eeq
where  $\bar{q}^{i,r} $ is the FWLLN limit of $\bar{Q}^{i,r}_n$ to be proved.
Let other scaled processes be defined similarly.  What we need are the FWLLN
$\bar{Q}^{i,r}_n  \Rightarrow \bar{q}^{i,r} $
and the associated FCLT
$\hat{Q}^{i,r}_n  \Rightarrow \hat{Q}^{i,r}$ in $D$ as $n \ra \infty$,
jointly with the limits in Theorem \ref{QFCLT}.  The extension to joint convergence with the other processes will be immediate
if the stochastic processes associated with new arrivals are independent of the initial conditions.
Otherwise, we require that we have the joint convergence
$(\hat{A}_n, \hat{Q}^{i,r}_n) \Rightarrow (\hat{A}, \hat{Q}^{i,r})$ in $D \times D$,
with the service times of new arrivals coming from a sequence of i.i.d. random variables, which is independent of
both the arrival processes and the initial conditions. We now give sufficient conditions to get these limits.

\paragraph{Assumptions for the Initial Conditions.}

\paragraph{Assumption 3:  i.i.d. service times.}  The service times of customers initially
in the system come from a sequence $\{\eta^i_j: j \ge 1\}$
of i.i.d. nonnegative random variables with a \textit{general} c.d.f. $F_i$ and $F_i(0) = 0$, independent of $n$ and
independent of the total number of customers initially present and
all random quantities associated with new arrivals.~~~\bsq

\paragraph{Assumption 4: independence and CLT for the initial number.}  The initial total number of customers in the system, $Q^{i,t}_n$,
is independent of the service times of the initial customers and all random quantities associated with new arrivals.
There exist (i) a nonnegative constant $\bar{q}^{i,t}$ and (ii) a random variable $\hat{Q}^{i,t}$ such that
\beql{f5}
\hat{Q}^{i,t}_n \equiv \frac{1}{\sqrt{n}}(Q^{i,t}_n -  n \bar{q}^{i,t})  \Rightarrow \hat{Q}^{i,t} \qinq \RR  \qasq n \ra \infty.~~~\bsq
\eeq

Paralleling Lemma \ref{repQn}, we have the representation result.

\begin{lemma} {\em $($representation of $Q^{i,r}_n$)}\label{repQin}
The process $Q^{i,r}_n$ can be represented as
\beql{E:Qin}
Q^{i,r}_n (y) = \sum_{j =1}^{Q^{i,t}_n}
\left(\mathbf{1}(\eta^{i}_{j} >y)-F_i^c (y) \right) +  Q^{i,t}_n F_i^c (y), \quad y\geq 0.
\eeq
\end{lemma}


\begin{theorem} {\em $($FWLLN and FCLT for the initial conditions$)$}\label{initialFCLT}
  Under Assumptions 3 and 4,
\begin{eqnarray}\label{f7}
\bar{Q}^{i,r}_n(y) & \Rightarrow & \bar{q}^{i,r}(y) \equiv  F_i^c (y) \bar{q}^{i,t}  \qinq D \qasq n \ra \infty,  \\
\hat{Q}^{i,r}_n (y)& \Rightarrow &\hat{Q}^{i,r}(y)\equiv  F_i^c (y) \hat{Q}^{i,t} + \sqrt{\bar{q}^{i,t}} B^0 (F_i(y)) \qinq D \qasq n \ra \infty, \nonumber
\end{eqnarray}
where $B^0$ is a Brownian bridge, independent of $\hat{Q}^{i,t}$.
\end{theorem}

We can combine Theorems \ref{QFWLLN}, \ref{QFCLT} and \ref{initialFCLT} to
 treat the total number of customers in the system at time $t$ with residual service times strictly greater than $y$,
 which we denote by $Q^{T,r}_n (t,y)$.
The key representation is
\beql{f9}
Q^{T,r}_n (t,y) = Q^r_n (t,y) + Q^{i,r} (t + y), \quad t \ge 0, \quad y \ge 0.
\eeq

\begin{coro}{\em $($FWLLN and FCLT for all customers$)$}\label{totalCor}
Under Assumptions $1-4$,
\begin{eqnarray}\label{f10}
\bar{Q}^{T,r}_n (t,y) & \equiv &  \bar{Q}^{i,r}_n (t + y) + \bar{Q}^r_n (t,y)   \Rightarrow  \bar{q}^{T,r}(t,y)
\equiv \bar{q}^{i,r} (t+y) + \bar{q}^r (t,y)   \nonumber \\
&& \quad \quad = F_i^c (t+y)  \bar{q}^{i,t} + \int_{0}^{t} F^c (t + y - s) \, d\bar{a} (s),  \\
\hat{Q}^{T,r}_n (t,y) & \equiv & \hat{Q}^{i,r}_n (t + y) + \hat{Q}^r_n (t,y)
\Rightarrow \hat{Q}^{T,r} (t,y) \equiv \hat{Q}^{i,r} (t+y) + \hat{Q}^r (t,y)   \nonumber \\
&& \quad \quad = F_i^c (t+y)  \hat{Q}^{i,t} + \sqrt{\bar{q}^{i,t}} B^0 (F_i (t+ y)) + \hat{X}^{c,r}_1 (t,y)  + \hat{X}^{c,r}_2 (t,y) + \hat{X}^{d,r} (t,y),  \nonumber
\end{eqnarray}
in $D_D$ as $n\ra\infty$, where $\hat{X}^{c,r}_1$, $\hat{X}^{c,r}_2$ and $ \hat{X}^{d,r}$ are given in {\em \eqn{E:X1c}}.
\end{coro}

\section{Proof of the FCLT}\label{secProof}

We now prove the FCLT in Theorem \ref{QFCLT}.
First, the joint convergence of
the processes
$$(\hat{A}_n, \hat{A}^c_n, \hat{A}^d_n, \{ \hat{A}^d_{n,i}:i\geq 1\}) \Rightarrow (\hat{A}, \hat{A}^c, \hat{A}^d, \{ \hat{A}^d_{i}:i\geq 1\})$$
follows from Theorem 9.5.1 in \cite{W02}. For the subsystem
with discrete service-time distribution, the limits follow from
an easy extension of \cite{GW91}. In \cite{GW91}, the
convergence to the limit $\hat{X}^{d,r}(t,y)$ is proved in the
space $D$ for each fixed $y\geq 0$, however, the convergence
can be easily generalized to be in the space $D_D$ since the
limit process $\hat{A}$ is assumed to be continuous here
(Assumption 1). Since the prelimit process of $\hat{X}^{d,r}$ is 
$$
\hat{X}_n^{d,r}(t,y) = \sum_{i=1}^{\infty} (\hat{A}^d_{n,i}(t) - \hat{A}^d_{n,i} (t - (\bar{x}_i - y)^{+}) ), \quad t, y \geq 0,
$$
it suffices to show that the mapping $\phi: D\ra D_D$ defined
by
$$
\phi(z)(t,y) \equiv \sum_{i=1}^{\infty} (z(t) - z(t-(\bar{x}_i - y)^{+})
$$
is continuous in the Skorohod $J_1$ topology and then apply the
continuous mapping theorem.
Moreover, in order to prove
$\hat{W}^{r,d}_n(t,y) \Rightarrow \hat{W}^{r,d}(t,y) $ in
$D_D$, where $\hat{W}^{r,d}_n(t,y)$ can be written as
\begin{eqnarray*}
\hat{W}^{r,d}_n(t,y) = \sum_{i=1}^{\infty} \int_y^{\bar{x}_i} (\hat{A}^d_{n,i}(t)
- \hat{A}^d_{n,i} (t - (\bar{x}_i - x)^{+}) ) d x , \quad t, y\geq 0,
\end{eqnarray*}
we need to prove the continuity of the mapping $\psi:D\ra D_D$ defined by
\begin{eqnarray*}
\psi(z)(t,y) = \int_y^{\bar{x}_i} (z(t) - z(t - (\bar{x}_i  - x)^{+})) d x, \quad z \in D, \quad t, y \geq 0.
\end{eqnarray*}
Since the limit $\hat{A}$ is continuous, it suffices to show the
uniform continuity of the mapping $\psi$ on compact intervals,
which follows from a direct argument.
Thus, we will only focus on the subsystem with continuous service-time distributions.
For notational convenience, we will simply suppose that $F$ in Assumption 2 is continuous such that $F_c = F$,
$\hat{A}_n^c = \hat{A}_n$ and similarly for other processes. In particular, we write $\hat{X}^{c,r}_1$
and $\hat{X}^{c,r}_2$ simply as $\hat{X}_1$ and $\hat{X}_2$, respectively.

One might hope to obtain a very fast proof by applying the continuous mapping theorem with an
appropriate continuous mapping.
That would seem to be possible, because
both the initial stochastic integral in \eqn{E:QnK} and the representation in Lemma \ref{repQn}
show that the scaled residual service queue-length process $\hat{Q}_n^r$ can be
regarded as the image of a deterministic function $h:D \times D_D \ra D_D$ mapping $(\hat{A}_n, \hat{K}_n)$
into $\hat{Q}_n^r$.
Given that $(\hat{A}_n, \hat{K}_n) \Rightarrow (\hat{A}, \hat{K})$ under Assumptions 1 and 2, we
would expect that corresponding limits for $\hat{Q}_n^r$ and the other processes would follow directly from an appropriate
continuous mapping theorem.  Unfortunately, the connecting map is complicated, being in the form
of a stochastic integral, with the limit of the component $\hat{X}_{n,2}$ involving a two-dimensional stochastic integral.
In fact, we will show below that we can easily treat the component $\hat{X}_{n,1}$ via the representation \eqn{Xn1parts}.
However, $\hat{X}_{n,2}$ presents a problem.
Unfortunately, the general results of weak convergence of stochastic integrals and
differential equations in \cite{KP91,M86,KP96} does not seem to apply.
Thus, instead, we will follow \cite{KP97} and prove the convergence in the classical way, by proving
tightness and convergence of the finite-dimensional distributions. (See \cite{RT09} for a different way.)

For us, the first step is to get convergence for the process $\hat{R}_n$ jointly with $(\hat{A}_n, \hat{K}_n)$ by exploiting
the composition map for a random time change, paralleling \S 13.2 of \cite{W02}; see \cite{TW08} for extensions to $D_D$.
Starting from $(\hat{A}_n, \hat{K}_n) \Rightarrow (\hat{A}, \hat{K})$,
we first obtain $(\hat{A}_n, \bar{A}_n, \hat{K}_n) \Rightarrow (\hat{A}, \bar{a}, \hat{K})$ by
applying \eqn{b2} and Theorem 11.4.5 of \cite{W02}.  We then apply the continuous mapping theorem for
composition applied in the space $D_D$, where the composition is with respect to the first component of $\hat{K}_n$,
and the limit $\bar{a}$ and $\hat{K}$ are both continuous (in the first component for $\hat{K}$).  That yields
\beql{h1}
(\hat{A}_n, \bar{A}_n, \hat{K}_n, \hat{R}_n) \Rightarrow (\hat{A}, \bar{a}, \hat{K}, \hat{R}) \qinq D^2 \times D_D^2,
\eeq
where $\hat{R}(t,x) = \hat{K}(\bar{a}(t),x) = U(\bar{a}(t), F(x))$ for $t\geq 0$ and $x\geq 0$.
Since $\hat{R}$ does not involve $\hat{A}$, we see that $\hat{A}_n$ and $\hat{R}_n$ are asymptotically independent.
Necessarily, then the processes $\hat{X}_{n,1}$ and $\hat{X}_{n,2}$ are asymptotically independent as well.

We use the classical method for establishing the limit
\beql{h2}
(\hat{A}_n, \bar{A}_n, \hat{K}_n, \hat{R}_n, \hat{X}_{n,1}, \hat{X}_{n,2})
\Rightarrow (\hat{A}, \bar{a}, \hat{K}, \hat{R}, \hat{X}_{1}, \hat{X}_{2})
\eeq
in $D^2 \times D_D^4$:  We show convergence of the finite-dimensional distributions
and tightness.  We get tightness for $\{(\hat{A}_n, \bar{A}_n, \hat{K}_n, \hat{R}_n): n \ge 1\}$ from the convergence in \eqn{h1}.
We use the fact that tightness on product spaces is equivalent to tightness on each of
the component spaces; see Theorem 11.6.7 of \cite{W02}.  Since we can write $\hat{X}_{n,1}$ as \eqn{Xn1parts},
the tightness and convergence of $\hat{X}_{n,1} \Rightarrow \hat{X}_1$ in $D_D$  can be obtained directly by
applying continuous mapping theorem if we can prove the mapping defined in \eqn{Xn1parts}
from $\hat{A}_n$ to $\hat{X}_{n,1}$ is continuous in $D_D$.  We will prove the continuity of this mapping in $D_D$
in \S \ref{secContinuity}.  We then establish tightness for $\{(\hat{X}_{n,1}, \hat{X}_{n,2}): n\geq 1\}$ in \S \ref{secTight}
and  the required convergence of the finite-dimensional distributions associated
with $\{(\hat{X}_{n,1}, \hat{X}_{n,2}): n\geq 1\}$  in \S \ref{secFiDi}.
Given the limit in \eqn{h2}, the rest of the limits in parts $(a)$ and $(b)$ follows from the continuous mapping theorem.
The limit in part (c) is an application of convergence preservation for composition with linear centering
as in Corollary 13.3.2 of \cite{W02}.  The component limits require finite second moments.

\subsection{Continuity of the Representation for ${\hat{X}_{n,1}}$ in $D_D$}\label{secContinuity}

In this section, we prove the continuity of the mapping $\phi: D \rightarrow  D_D$ defined by
\begin{eqnarray}  \label{E:Xn1map}
\phi(x)(t,y) &\equiv&  F^c (y)  x(t)- \int_0^t  x(s-)d F(t + y -s),
\end{eqnarray}
for $x \in D$ and $t,y\geq 0$.
By \eqn{Xn1parts} and \eqn{E:X1c}, we have $\hat{X}_{n,1}(t,y) = \phi(\hat{A}_n)(t,y)$ and $\hat{X}_1(t,y) = \phi(\hat{A})(t,y)$.

\begin{lemma} \label{Xn1mapcont}
The mapping $\phi$ defined in $\eqn{E:Xn1map}$ is continuous in $D_D$.
\end{lemma}

\paragraph{Proof.}

Suppose that $x_n \rightarrow x$ in $D$.  We need to show that
$d_{D_D} (\phi(x_n), \phi(x))\rightarrow 0$ as $n \rightarrow \infty$.
Let $T>0$ be a continuity point of $x$ and consider the time domain $[0,T] \times [0,\infty)$.
By the convergence $x_n\ra x$ in $(D,J_1)$ as $n\ra\infty$, there exist increasing homeomorphisms $\lambda_n$
of the interval $[0,T]$ such that $|| x_n - x\circ \lambda_n ||_T \ra 0$ and $|| \lambda_n - e||_T \ra 0$ as $n\ra\infty$,
where $e(t) = t$ for all $t\geq 0$ and $||y ||_T = \sup_{t\in[0,T]} |y(t)|$ for any $y \in D$.
Let $M = \sup_{0\leq t \leq T} |x(t)| < \infty$.
Since $F$ is continuous, it suffices to show that
\begin{eqnarray*}
&& ||\phi(x_n) ( \cdot, \cdot) -  \phi(x)(\lambda_{n}(\cdot), \cdot) ||_{T}  \\
&=& \sup_{(t,y) \in [0,T]\times[0,\infty)} |\phi(x_n)(t,y) - \phi(x)(\lambda_{n}(t),y)|  \rightarrow 0, \qasq  n \rightarrow \infty.
\end{eqnarray*}
Now, we have
\begin{eqnarray*}
&&|\phi(x_n)(t, y) - \phi(x)(\lambda_{n}(t),y)|  \\
&=& \Big|  F^c(y) x_n(t) - \int_0^{t}  x_n(s-)d F(t +y-s) {}\\
&& {} \quad -  F^c(y) x(\lambda_{n}(t))  + \int_0^{\lambda_{n}(t)}  x(s-)d F(\lambda_{n}(t) + y -s) \Big|\\
&\leq&   F^c(y)\big| x_n(t) -x(\lambda_{n}(t))\big|  {}\\
&&{}  + \Big|  \int_0^{t}  x_n(s-)d F(t +y-s)
-  \int_0^{\lambda_{n}(t)}  x(s-)d F(\lambda_{n}(t) + y -s)\Big| \\
&=&  F^c(y)\big| x_n(t) -x(\lambda_{n}(t))\big| {} \\
&& +  \Big|  \int_0^{t}  x_n(s-) d F(t +y-s)   -  \int_0^t  x(\lambda_{n}(s)-)d  F(\lambda_{n}(t) + y -\lambda_{n}(s)) \Big| \\
&\leq& F^c(y)\big| x_n(t) -x(\lambda_{n}(t))\big| +
 \Big|\int_0^t ( x_n(s-) - x(\lambda_n(s)-) )   d F(t + y -s)\Big|  {}\\
&&{} + \Big| \int_0^t x(\lambda_n(s)-) d (F(\lambda_{n}(t) +y-\lambda_n(s)) - F(t + y -s))  \Big| {}\\
&\leq& F^c(y)\big| x_n(t) -x(\lambda_{n}(t))\big|  + ||x_n - x \circ \lambda_n ||_T |F(y) - F(t+y)|\\
&&{} + M | F(\lambda_n(t) +y) - F(t+y)|  {}\\
&\leq&  3||x_n - x \circ \lambda_n ||_T  + M | F(\lambda_n(t) +y) - F(t+y)| .
\end{eqnarray*}
The third term in the third inequality follows from the uniform continuity of the integrator because $F$ is continuous, monotone and bounded.
By taking the supremum over $(t,y) \in [0,T]\times[0,\infty)$, the first term converges to $0$ by the convergence of $x_n\ra x$ in $D$, and
the second term converges to $0$ by the uniform convergence of $\lambda_n \ra e$ in $[0,T]$ and the continuity of $F$.
This implies the initial convergence to be shown, so that the mapping $\phi: D\rightarrow D_D$ is indeed continuous. ~~~\bsq

\subsection{Tightness }\label{secTight}

In this section, we establish tightness for the sequence of scaled processes in \eqn{E:fclt}.
It suffices to prove tightness of the sequences of processes
 $\{\hat{X}_{n,1}: n\geq1\}$ and  $\{\hat{X}_{n,2}: n\geq1\}$ in $D_D$. By Assumption 1,
 the sequence of processes $\{\hat{A}_n:n\geq 1\}$ is tight.  The tightness of $\{\hat{X}_{n,1}\}$
 follows from the continuity of the mapping $\phi$ in $D_D$.  It remains to show the tightness of $\{\hat{X}_{n,2}\}$
 and then we obtain tightness of the sequences of processes $\{\bar{Q}_n^r:n\geq1\}$ and $\{\hat{Q}_n^r:n\geq1\}$
 using the fact that tightness of product spaces is equivalent to the tightness on each of the component spaces.

\begin{theorem} \label{TQ}
 Under Assumptions 1 and 2 $($F is continuous$)$, the sequence of processes $\{\hat{X}_{n,1}:n\geq 1\}$, $\{\hat{X}_{n,2}:n\geq 1\}$,
 $\{\bar{Q}_n^r:n\geq1\}$ and $\{\hat{Q}_n^r:n\geq1\}$ are individually and jointly tight.
\end{theorem}

In order to prove the tightness of $\{\hat{X}_{n,2}: n\geq 1\}$ defined in \eqn{E:Mn2},
we will closely follow the approach in \cite{KP97}, but we must adjust to the tightness criteria in $D_D$.
The following tightness criteria are adapted to $D_D$ from Theorem 3.8.6 of Ethier and Kurtz \cite{EK86}.
For a review of tightness criteria for processes in the space $D$, see \cite{WW07}.

\begin{theorem} \label{TCDD}

A sequence of stochastic processes $\{X_n: n\geq 1\}$ in $D_D$ is tight if and only if

(i) the sequence  $\{X_n: n \geq 1\}$  is stochastically bounded in $D_D$, i.e.,
for all $\epsilon > 0$, there exists a compact subset $K\subset \mathbb{R}$ such that
$$
P(||X_n||_T \in K) > 1-\epsilon, \quad \rm{for\quad all} \quad n \geq 1,
$$
where $||X_n||_T = \sup_{s\in[0,T]}\{\sup_{t \in [0,T]} |X_n(s,t)| \}$;

and any one of the following

(ii)
For all $\delta > 0$, and all uniformly bounded sequences   $\{\tau_n: n\geq 1\}$ where
for each $n$, $\tau_n$ is a stopping time  with respect to the natural
filtration $\mathbf{F}_n = \{\mathcal{F}_n(t), t\in[0,T]\}$ where $\mathcal{F}_n(t)
= \sigma\{X_n(s,\cdot): 0 \leq s \leq t\}$, there exists a constant $\beta >0$ such that
\beq
\lim_{\delta\downarrow 0}\limsup_{n\rightarrow\infty} \sup_{\tau_n}E[(1\wedge d_{J_1}(X_n(\tau_n+\delta, \cdot), X_n(\tau_n,\cdot)))^{\beta}] = 0;
\eeqno

or

(ii')
For all $\delta > 0$, there exist a constant $\beta$ and random variables $\gamma_n(\delta) \geq 0$ such that for each $n$, w.p.1,
$$
E[(1\wedge d_{J_1}(X_n(s+u, \cdot),  X_n(s,\cdot)))^{\beta}| \mathcal{F}_{n}](1\wedge  d_{J_1} ( X_n(s-v, \cdot),  X_n(s,\cdot)))^{\beta}\leq E[\gamma_n(\delta) | \mathcal{F}_{n}],
$$
for all $0 \leq s \leq T$, $0\leq u \leq \delta$ and $0 \leq v \leq s \wedge \delta$,
where  $\mathbf{F}_n =\{ \mathcal{F}_{n}(t):  t \in [0,T] \}$ with  $\mathcal{F}_n(t) = \sigma\{X_n(s,\cdot): 0 \leq s \leq t\}$
and
\beq
\lim_{\delta\downarrow 0} \limsup_{n\rightarrow\infty} E[\gamma_n(\delta)] = 0.
\eeqno

\end{theorem}

\begin{remark} The following condition is sufficient, but not necessary, for condition $(ii)$ in Theorem \ref{TCDD}:

\textit{
For all $\delta_n \downarrow 0$ and for all uniformly bounded sequences  $\{\tau_n: n\geq 1\}$, where
for each $n$, $\tau_n$ is a stopping time  with respect to the natural
filtration $\mathbf{F}_n =\{ \mathcal{F}_{n}(t):  t \in [0,T] \}$ with  $\mathcal{F}_n(t) = \sigma\{X_n(s,\cdot): 0 \leq s \leq t\}$,
\beq
d_{J_1}(X_n(\tau_n + \delta_n, \cdot),  X_n(\tau_n, \cdot) ) \Rightarrow 0, \qasq n\ra \infty. 
\eeqno
}
\end{remark}

We will also need to generalize the tightness criteria in Lemma VI.3.32 in \cite{JS87} for processes in the space $D$ to those in the space $D_D$ as in the following lemma, and its proof also follows from that in \cite{JS87} with inequalities for the modulus of continuity for functions in the space $D_D$.

\begin{lemma} \label{DDtightJS}
Suppose that a sequence of processes $\{X_n:n\geq 1\}$ in the space $D_D$ can be decomposed into two sequences $\{Y^{q}_n:n\geq 1\}$ and $\{Z^q_n:n\geq 1\}$ for some parameter $q\in \mathbb{N}$, i.e., $X_n = Y^q_n + Z^q_n$ for each $n\geq1$, and that $(i)$ the sequence $\{Y^{q}_n:n\geq 1\}$ is tight in the space $D_D$ and $(ii)$ for all $T>0$ and $\delta>0$, $\lim_{q\ra\infty} \limsup_{n\ra\infty} P(\sup_{t,y\leq T} |Z^q_n(t,y)| > \delta) =0$.  Then, the sequence $\{X_n:n\geq 1\}$ is tight in the space $D_D$.
\end{lemma}

We first give a decomposition of the process $\hat{X}_{n,2}$ for each $n$.
Following \cite{KP97}, we can write $\hat{R}_n(t,y)$ in \eqn{E:Rn} as
\beq
\hat{R}_n(t,y) = -\int_0^y \frac{\hat{R}_{n}(t,x)}{1-F(x)} d F(x) + \hat{L}_n(t,y),
\eeqno
where
\beq
\hat{L}_n(t,y) = \frac{1}{\sqrt{n}}\sum\limits_{i =1}^{A_n(t)}\Big(\mathbf{1}(\eta_i \leq y)-\int_0^{y\wedge\eta_i}
\frac{1}{1-F(x)} d F(x)\Big).
\eeqno
We remark that we need not consider the left-hand limit of $\hat{R}_n$ in the second argument, as was done in \cite{KP97},
since the service-time c.d.f  $F$ is assumed to be continuous, while $F$ is allowed to be discontinuous in \cite{KP97}.
Hence, $\hat{X}_{n,2}$ can be written as
\begin{equation}
\hat{X}_{n,2}(t,y) = \hat{G}_n(t,y) + \hat{H}_n(t,y), \qforq t\geq 0 \qandq y\geq 0,
\end{equation}
where
\begin{eqnarray}\label{E:hGnhRn}
\hat{G}_n(t,y) &\equiv& \int_0^t \int_0^{\infty} \mathbf{1}( s+x \leq t+y) d\Big(
-\int_0^x \frac{\hat{R}_{n}(s,v)}{1-F(v)} d F(v) \Big) \nonumber \\
&=&  - \int_0^{t+y} \frac{\hat{R}_{n}(t+y-x,x)}{1-F(x)} d F(x),
\end{eqnarray}
and
\begin{equation}
\hat{H}_n(t,y) \equiv \int_0^t \int_0^{\infty}\mathbf{1}( s+x \leq t+y)d \hat{L}_n(s,x).
\end{equation}
Thus, the tightness of $\{\hat{X}_{n,2}\}$ follows from the tightness of  $\{\hat{G}_n\} $ and $\{\hat{H}_n\} $.
We will establish their tightness in the following two lemmas.

\begin{lemma} \label{TGn}
 Under Assumptions 1 and 2 $($F is continuous$)$, the sequence of processes
 $\{\hat{G}_n: n\geq 1\} \equiv \{ \{\hat{G}_n(t,y): t\geq 0, y\geq 0\}, n\geq 1\}$ is tight in the space $D_D$.
\end{lemma}

\paragraph{Proof.}
We will apply Lemma \ref{DDtightJS}.   We define the sequence of processes $\{\hat{G}_n^{\epsilon}: n\geq 1\}$, for some $\epsilon \in (0,1)$, by
\begin{equation}\label{E:TGne1}
\hat{G}_n^{\epsilon}(t,y) \equiv
- \int_0^{t+y} \frac{\hat{R}_{n}(t+y-x,x)}{1-F(x)} \mathbf{1}(F(x) \leq 1-\epsilon) d F(x),\quad t,y \geq 0.
\end{equation}
We will prove that  $\{\hat{G}_n^{\epsilon}: n\geq 1\}$ is tight in $D_D$ and
\begin{equation}\label{E:TGnebb}
\lim_{\epsilon\downarrow 0} \limsup_n P\Big( \sup_{t,y \leq T} \Big |
 \int_0^{t+y} \frac{\hat{R}_{n}(t+y-x,x)}{1-F(x)} \mathbf{1}(F(x) > 1-\epsilon) d F(x)\Big | > \delta \Big) =0,
\end{equation}
for each $\delta >0$ and $T>0$,
and thus will conclude that the sequence $\{\hat{G}_n\}$ is tight in $D_D$ by Lemma \ref{DDtightJS}.
 It is easy to see that \eqn{E:TGnebb} follows easily from (3.23) in \cite{KP97}. So we only need to prove the tightness of the sequence of processes
$\{\hat{G}_n^{\epsilon}: n\geq 1\}$.

%

Recall that $\hat{R}_n(t+y-x,x) = \hat{U}_n(\bar{A}_n(t+y-x), F(x))$.  By \eqn{b2}  and $\hat{U}_n\Rightarrow U$ in \eqn{b5} as $n\ra\infty$,
and by applying the continuous mapping theorem to the composition map of $\hat{U}_n$ with respect
to the first argument (Theorem 13.2.2, \cite{W02}), we obtain
$$
\hat{R}_n(t+y-x, x) =  \hat{U}_n(\bar{A}_n(t+y-x), F(x)) \Rightarrow U(\bar{a}(t+y-x),F( x)) \qinq D_D,
$$
as $n\ra\infty$.
The weak convergence of $\{\hat{R}_n:n\geq 1\}$ in $D_D$ implies that  $\{\hat{R}_n:n\geq 1\}$ is stochastically bounded,
so the integral representation of $\hat{G}^{\epsilon}_n$ in terms of $\hat{R}_n$ in \eqn{E:TGne1} implies that
  $\{\hat{G}^{\epsilon}_n: n\geq1\}$ is also stochastically bounded in $D_D$. We apply Theorem \ref{TCDD}
  to prove the tightness of $\{\hat{G}_n^{\epsilon}: n\geq 1\}$ in $D_D$. In this case, it is convenient
  to use the sufficient criterion in the remark  right after Theorem \ref{TCDD}.

Let $\mathbf{G}_n = \{\mathcal{G}_n(t): t \in [0,T] \}$ be a filtration defined by
\begin{eqnarray*}
\mathcal{G}_n(t)& =& \sigma\{\hat{R}_n(s,\cdot): 0 \leq s \leq t\} \vee \mathcal{N} \\
&=& \sigma\{\eta_i \leq  x: 1\leq i \leq A_n(t), x \geq 0\} \vee \sigma\{A_n(s): 0 \leq s \leq t\}\vee \mathcal{N},
\end{eqnarray*}
where $\mathcal{N}$ includes all the null sets. Note that the filtration  $\mathbf{G}_n$ satisfies the usual conditions (Chapter 1, \cite{KS91} and proof of Lemma 3.1 in \cite{KP97}).
Let $\delta_n \downarrow 0$ and   $\{\tau_n: n\geq 1\}$ be a uniformly bounded sequence,
where for each $n$, $\tau_n$ is a stopping times  with respect to the filtration $\mathbf{G}_n$.
Then, it suffices to show that
\beq
d_{J_1} (\hat{G}_n^{\epsilon}(\tau_n+ \delta_n, \cdot),  \hat{G}_n^{\epsilon}(\tau_n, \cdot)) \Rightarrow 0, \qasq n\rightarrow \infty.
\eeqno
Consider  any sequence of nondecreasing homeomorphism $\{\lambda_n: n\geq 1\}$ on $[0,T]$
such that $\lim_{n\rightarrow \infty} \lambda_n(y) = y$ uniformly in $y \in[0,T]$.
We want to show that the following holds:
\beq
\sup_{0 \leq y \leq T} \Big| \hat{G}_n^{\epsilon}(\tau_n+ \delta_n, \lambda_n(y)) - \hat{G}_n^{\epsilon}(\tau_n, y) \Big|
\Rightarrow 0, \qasq n \ra\infty.
\eeqno
Now,
\begin{eqnarray*}
&&\sup_{0 \leq y \leq T} \left| \hat{G}_n^{\epsilon}(\tau_n+ \delta_n, \lambda_n(y))
- \hat{G}_n^{\epsilon}(\tau_n, y) \right| \\
&=& \sup_{0 \leq y \leq T} \Big|   \int_0^{\tau_n+ \delta_n+\lambda_n(y)} \frac{\hat{R}_{n}(\tau_n+ \delta_n+\lambda_n(y)-x,x)}{1-F(x)} \mathbf{1}(F(x)
\leq 1-\epsilon) d F(x)   \nonumber \\
& & {}  -   \int_0^{\tau_n+y} \frac{\hat{R}_{n}(\tau_n+y-x,x)}{1-F(x)} \mathbf{1}(F(x) \leq 1-\epsilon) d F(x)  \Big|\\
&\leq&   \sup_{0 \leq y \leq T} \Big|   \int_0^{\tau_n+ \delta_n+\lambda_n(y)} \frac{\hat{R}_{n}(\tau_n+ \delta_n+\lambda_n(y)-x,x)
-\hat{R}_{n}(\tau_n+y-x,x)}{1-F(x)} \mathbf{1}(F(x) \leq 1-\epsilon) d F(x) \Big| {}\\
&&{}    +   \sup_{0 \leq y \leq T} \Big|   \int_{0}^{\tau_n+ \delta_n+\lambda_n(y)} \frac{\hat{R}_{n} (\tau_n+y-x,x)}{1-F(x)}\mathbf{1}(F(x) \leq 1-\epsilon) d F(x) {}\\
&& {} - \int_0^{\tau_n + y}   \frac{\hat{R}_{n} (\tau_n+y-x,x)}{1-F(x)}\mathbf{1}(F(x) \leq 1-\epsilon) d F(x) \Big|\\
&\Rightarrow& 0
\end{eqnarray*}
as $n\rightarrow \infty$, where the first and the second terms converge to 0 by
the stochastic boundedness and weak convergence of $\hat{R}_n$ in $D_D$,
and because $\tau_n$ is uniformly bounded, $\lambda_n(y)$ converges to $y$ uniformly in $[0,T]$, and $\delta_n \downarrow 0$ as $n\ra\infty$.
Hence, the processes $\{\hat{G}_n^{\epsilon}\}$ are tight in $D_D$ and the proof is completed.
~~~\bsq

\begin{lemma}\label{THn}
 Under Assumptions 1 and 2 $($F is continuous$)$, the sequence of processes $\{\hat{H}_n: n\geq 1\} \equiv \{ \{\hat{H}_n(t,y): t\geq 0, y\geq 0 \}, n\geq 1\}$
 is tight in $D_D$.

\end{lemma}

\paragraph{Proof.}
As in Lemma 3.7 in \cite{KP97}, we write the process $\hat{H}_n$ as
\beq
\hat{H}_n(t,y) =  \frac{1}{\sqrt{n}} \sum_{i = 1}^{A_n(t)} \left(\mathbf{1}(0 \leq \eta_i \leq t+y- \tau_i^n)
- \int_{0}^{\eta_i \wedge (t+y-\tau_i^n)^{+}} \frac{1}{1-F(u)} d F(u)\right).
\eeqno
We will apply Theorem \ref{TCDD} to prove the tightness of $\{\hat{H}_n: n\geq 1\}$ in $D_D$. In this case,
it is convenient to use criterion $(ii)$ in Theorem \ref{TCDD}.  We will first prove that  this criterion holds,
and then prove the stochastic boundedness of the sequence of processes $\{\hat{H}_n: n \geq 1\}$.

Let $\mathbf{H}_n = \{\mathcal{H}_n(t): t \in [0,T]\}$  be a filtration defined by
\begin{eqnarray*}
\mathcal{H}_n(t)& =& \sigma\{\hat{H}_n(s,\cdot): 0 \leq s \leq t\}\vee \mathcal{N}  \\
&=& \sigma\{\eta_i \leq s+x-\tau_i^n: 1\leq i \leq A_n(t),  x \geq 0, 0 \leq s \leq t \} \vee \{A_n(s): 0 \leq s \leq t\}\vee \mathcal{N},
\end{eqnarray*}
where $\mathcal{N}$ includes all the null sets.  The filtration $\mathbf{H}_n$ satisfies the usual conditions (see p. 254 in \cite{KP97}).

Let $\delta > 0$ and  $\{\kappa_n: n\geq 1\}$  be a uniformly bounded sequence, where  for each $n$, $\kappa_n$ is a stopping time
with respect to the  filtration $\mathbf{H}_n$.  It suffices to show that
\begin{equation} \label{E:THne1}
\lim_{\delta\downarrow 0}\limsup_{n\rightarrow\infty} \sup_{\kappa_n}E[d_{J_1}( \hat{H}_n(\kappa_n+\delta, \cdot), \hat{H}_n(\kappa_n,\cdot)) ^2] = 0.
\end{equation}

Consider any sequence of nondecreasing homeomorphism $\{\lambda_n: n\geq 1\}$ on $[0,T]$
such that $\lim_{n\rightarrow \infty} \lambda_n(y) = y$ uniformly in $y \in[0,T]$.
We want to show that the following holds:
\begin{equation} \label{E:THne2}
\lim_{\delta\downarrow 0}\limsup_{n\rightarrow\infty} \sup_{\kappa_n}E\Big[ \Big( \sup_{0\leq y \leq T}| \hat{H}_n(\kappa_n+\delta, \lambda_n(y))
- \hat{H}_n(\kappa_n,y)| \Big)^2 \Big] = 0.
\end{equation}

Define the processes $\hat{H}_{n,i} \equiv \{\hat{H}_{n,i}(t,y): t,y \geq 0\}$ by
\beq
\hat{H}_{n,i}(t,y) \equiv \mathbf{1}(0 \leq \eta_i \leq t+y- \tau_i^n)
- \int_{0}^{\eta_i \wedge (t+y-\tau_i^n)^{+}} \frac{1}{1-F(u)} d F(u).
\eeqno

As in Lemma 3.5 in \cite{KP97}, one can check that for each fixed $y$ and for each $i$,
the process $\{\hat{H}_{n,i}(t,y):t \geq 0\}$ is a square integrable martingale
with respect to the filtration $\mathbf{H}_n $ and it has predictable quadratic variation
\beq
\langle \hat{H}_{n,i} (\cdot, y)\rangle (t)
 = \langle \hat{H}_{n,i} \rangle (t,y) = \int_{0}^{\eta_i \wedge (t+y-\tau_i^n)^{+}} \frac{1}{1-F(u)} d F(u), \qforq t\geq 0,
\eeqno
and that the $\mathbf{H}_n$ martingales $\hat{H}_{n,i}(\cdot, y)$
and $\hat{H}_{n,j}(\cdot, y)$ for each fixed $y$ are orthogonal for $i \neq j$.

Thus, for each fixed $y$ and constant $K>0$,  the process  $\hat{H}_{n}^{(K)} = \{\hat{H}_{n}^{(K)}(t, y): t\geq 0\}$ defined by
\begin{eqnarray}
\hat{H}_{n}^{(K)}(t,y) & = &  \frac{1}{\sqrt{n}} \sum_{i = 1}^{n ( \bar{A}_n(t)\wedge  K) } \Big(\mathbf{1}(0 \leq \eta_i \leq t+y- \tau_i^n)  - \int_{0}^{\eta_i \wedge (t+y-\tau_i^n)^{+}} \frac{1}{1-F(u)} d F(u)\Big), \nonumber
\end{eqnarray}
 is an $\mathbf{H}_n$ square integrable martingale with predictable quadratic variation
 \beq
\langle \hat{H}_{n}^{(K)} (\cdot, y)\rangle (t)  = \langle \hat{H}_{n}^{(K)} \rangle (t,y)
= \frac{1}{n}\sum_{i = 1}^{n (\bar{A}_n(t)\wedge K) }  \int_{0}^{\eta_i \wedge (t+y-\tau_i^n)^{+}} \frac{1}{1-F(u)} d F(u),
\eeqno
for $t\geq 0$.
By the SLLN,
\begin{equation}\label{E:Thne3}
\frac{1}{n} \sum_{i = 1}^{\lfloor nt \rfloor}  \int_{0}^{\eta_i} \frac{1}{1-F(u)} d F(u)
\rightarrow t, \quad \it{a.s.}  \qasq  n \ra \infty.
\end{equation}
So for each fixed $y$, the sequence of quadratic variations $\{\langle \hat{H}_{n}^{(K)} (\cdot, y)\rangle: n\geq 1\}$
is $C$-tight by the continuity of $\bar{a}$ (Recall that a sequence $\{Y_n\}$ is said to be C-tight if it is tight and the limit of any convergent subsequence must have continuous sample paths.).  It follows by Theorem 3.6 in \cite{WW07}
that the sequence $\{\hat{H}_{n}^{(K)}(\cdot, y):  n\geq 1\}$ is $C$-tight for each fixed $y$.

Now,  to prove \eqn{E:THne2}, we have
\begin{eqnarray*}
&& E\Big[ \Big( \sup_{0\leq y \leq T} \big| \hat{H}_n(\kappa_n+\delta, \lambda_n(y)) - \hat{H}_n(\kappa_n,y) \big| \Big)^2\Big] \\
&\leq& 2 E\Big[ \sup_{0\leq y \leq T}\big| \hat{H}_n(\kappa_n+\delta, \lambda_n(y))
- \hat{H}_n(\kappa_n, \lambda_n(y))\big|^2\Big]  {}\\
&&{} + 2E\Big[\sup_{0\leq y \leq T}\big| \hat{H}_n(\kappa_n,\lambda_n(y) )
- \hat{H}_n(\kappa_n,y)\big|^2\Big]  \\
 &=& 2 \lim_{K\ra\infty} E\Big[ \sup_{0\leq y \leq T}\big| \hat{H}_n^{(K)}(\kappa_n+\delta, \lambda_n(y))
- \hat{H}_n^{(K)}(\kappa_n,\lambda_n(y))\big|^2\Big] \\
&&+ 2  \lim_{K\ra\infty}E\Big[\sup_{0\leq y \leq T}\big| \hat{H}_n^{(K)}(\kappa_n, \lambda_n(y)) - \hat{H}_n^{(K)}(\kappa_n,y)\big|^2\Big],
\end{eqnarray*}
where the equality holds by the dominated convergence and by stochastic boundedness of $A_n$.
The first term converges to 0 as $n\ra\infty$ and $\delta \downarrow 0$ by the assumptions
on $\kappa_n$ and $\lambda_n$ and $C$-tightness of $\{\hat{H}_{n}^{(K)}:  n\geq 1\}$.  We
conclude that the second term converges to 0 by observing
\begin{eqnarray*}
&&\hat{H}_n^{(K)}(\kappa_n, \lambda_n(y)) - \hat{H}_n^{(K)}(\kappa_n,y)\\
&=&  \frac{1}{\sqrt{n}} \sum_{i = 1}^{A_n(\kappa_n) \wedge K} \Big(\mathbf{1}(0 \leq \eta_i \leq \kappa_n+\lambda_n(y)- \tau_i^n) - \mathbf{1}(0 \leq \eta_i \leq \kappa_n+y- \tau_i^n) {}\\
&& {} -\Big( \int_{0}^{\eta_i \wedge (\kappa_n+ \lambda_n(y)-\tau_i^n)^{+}} \frac{1}{1-F(u)} d F(u)\
- \int_{0}^{\eta_i \wedge (\kappa_n+y-\tau_i^n)^{+}} \frac{1}{1-F(u)} d F(u) \Big) \Big).
\end{eqnarray*}
Thus we obtain \eqn{E:THne2}.

Now we prove the stochastic boundedness of $\{\hat{H}_n: n \geq 1\}$ in $D_D$.  We observe that
for each $n$, each sample path of the process $\hat{H}_n$ is  bounded by that of the process $\tilde{H}_n$ defined by
$$
\tilde{H}_n(t,y) =  \frac{1}{\sqrt{n}} \sum_{i = 1}^{A_n(t+y)} \left(\mathbf{1}(0 \leq \eta_i \leq t+y- \tau_i^n)
- \int_{0}^{\eta_i \wedge (t+y-\tau_i^n)^{+}} \frac{1}{1-F(u)} d F(u)\right).
$$
The stochastic boundedness of $\{\tilde{H}_n: n \geq 1\}$ in $D_D$ follows directly from the proof of Lemma 3.7 in \cite{KP97}.
Therefore, $\{\hat{H}_n: n \geq 1\}$ is stochastically bounded, so that tightness of $\{\hat{H}_n: n\geq 1\}$ in $D_D$ is proved.~~~\bsq

\subsection{Convergence of the Finite-Dimensional Distributions} \label{secFiDi}

In this section, we complete the proof of the convergence $(\hat{X}_{n,1}, \hat{X}_{n,2})
\Rightarrow (\hat{X}_1, \hat{X}_2)$ in $D_D \times D_D$ by proving that the
finite-dimensional distributions of $(\hat{X}_{n,1}, \hat{X}_{n,2})$ converge to those of $(\hat{X}_1, \hat{X}_2)$
since we have proved the tightness of $\{(\hat{X}_{n,1}, \hat{X}_{n,2}): n\geq 1\} $ in \S \ref{secTight}.
We will mostly have to deal with $\hat{X}_{n,2}$, since we have already shown convergence of $\hat{X}_{n,1}$.
Our argument for $\hat{X}_{n,2}$ will also enable us to establish joint convergence of the two finite-dimensional distributions.

\begin{lemma}\label{Mn2C}
 Under Assumptions 1 and 2 $($F is continuous$)$,  the finite-dimensional distributions of $(\hat{X}_{n,1}, \hat{X}_{n,2})$
 converge to those of  $(\hat{X}_1, \hat{X}_2)$ as $n\rightarrow \infty$.
\end{lemma}

\paragraph{Proof.}

First of all, we understand the integrals $\hat{X}_{n,2}$ in \eqn{E:Mn2}
and $\hat{X}_2$ ($\equiv \hat{X}^{c,r}_2$) in \eqn{E:X1c}  as mean-square integrals, so that they can be represented as
\beq
\hat{X}_{n,2}(t,y) =
\textrm{l.i.m.}_{k\rightarrow\infty}\hat{X}_{n,2,k}(t,y), \qandq
\hat{X}_2 (t,y) = \textrm{l.i.m.}_{k\rightarrow\infty} \hat{X}_{2,k}(t,y),
\eeqno
where $ \textrm{l.i.m.}$ means limit in mean square, that is,
$$
\lim_{k\ra\infty}E[(\hat{X}_{n,2}(t,y) -\hat{X}_{n,2,k}(t,y))^2 ] =0 \qandq \lim_{k\ra\infty}E[(\hat{X}_{2}(t,y) -\hat{X}_{2,k}(t,y))^2 ] =0,
$$
\begin{eqnarray*}
\hat{X}_{n,2,k}(t,y) &\equiv& - \int_0^t\int_0^{\infty}\mathbf{1}^y_{k,t}(s,x)d \hat{U}_n(\bar{A}_n(s), F(x)) \nonumber\\
&=& - \sum_{i=1}^k \Big[ \Delta_{\hat{U}_n}(\bar{A}_n(s_{i-1}^k),\bar{A}_n(s_i^k), 0,F(t+y-s_i^k))\Big],
\end{eqnarray*}
and
\begin{eqnarray*}
\hat{X}_{2,k} (t,y)&\equiv& - \int_0^t\int_0^{\infty} \mathbf{1}^y_{k,t}(s,x)d U(\bar{a}(s), F(x)) \\
&=& - \sum_{i=1}^k \Big[ \Delta_{U}(\bar{a}(s_{i-1}^k),\bar{a}(s_i^k), 0, F(t+y-s_i^k))\Big],
\end{eqnarray*}
where
 $\mathbf{1}^y_{k,t}$ is defined by
\beql{E:M2e1y}
\mathbf{1}^y_{k,t}(s,x)= \mathbf{1}(s=0)\mathbf{1}(x\leq t+y) +
\sum_{i = 1}^{k} \mathbf{1}(s\in(s_{i-1}^k,s_{i}^k ])\mathbf{1}(x
\leq t +y- s_i^k),
\eeq
with the points $0 = s_0^k < s_1^k < ... < s_k^k = t$ chosen so
that $\max_{1\leq i \leq k} |s_{i-1}^k-s_{i}^k|\rightarrow 0$ as
$k\rightarrow \infty$,
 and $\Delta_{\hat{U}_n}$ and $\Delta_{U}$ are defined as $\Delta_{\hat{K}}$ in \eqn{Kinc}.

 We prove the convergence of the finite-dimensional distributions of $\hat{X}_{n,2}$
 to those of $\hat{X}_2$ by taking advantage of the convergence  of $\hat{U}_n
 \Rightarrow U$ as $n\ra\infty$ in $D([0,\infty),D([0,1], \RR) )$ (see \eqn{b5}),
 for which we define another process $\{\tilde{X}_{n,2,k}(t,y): t, y \geq 0\}$ in $D_D$
 for each $n$ by replacing the $\bar{A}_n$ terms in $\Delta_{\hat{U}_n}$ of  $\hat{X}_{n,2,k}$ by $\bar{a}$ as follows,
\begin{eqnarray*}
\tilde{X}_{n,2,k}(t,y) &\equiv& - \int_0^t\int_0^{\infty}\mathbf{1}^y_{k,t}(s,x)d \hat{U}_n(\bar{a}(s), F(x)) \\
&=& - \sum_{i=1}^k \Big[ \Delta_{\hat{U}_n}(\bar{a}(s_{i-1}^k), \bar{a}(s_i^k), 0, F(t+y-s_i^k))\Big].\nonumber
\end{eqnarray*}
Hence, we easily obtain the convergence of the finite-dimensional distributions of $\tilde{X}_{n,2,k}$
to those of $\hat{X}_{2,k}$ as $n\ra\infty$, since $\bar{a}$ and $F$ are both continuous by Assumptions 1 and 2,
and the finite-dimensional distributions of $\hat{U}_n$ converge to those of $U$ as $n\ra\infty$ and $U$ is continuous.

Moreover, since $\hat{K}_n$ $( \hat{U}_n)$ and $A_n$ are independent by Assumptions 1 and 2,
$\tilde{X}_{n,2,k}$ and $\hat{X}_{n,1}$  are independent, and since the limit processes $\hat{X}_{2,k}$
and $\hat{X}_1$ are also independent, we obtain the joint convergence
of the finite-dimensional distributions of $(\hat{X}_{n,1}, \tilde{X}_{n,2,k})$  to those of $(\hat{X}_1, \hat{X}_{2,k})$ as $n\ra\infty$.

Now it suffices to show that the difference between $\hat{X}_{n,2,k}$ and $\tilde{X}_{n,2,k}$ is asymptotically negligible
in probability as $n\ra\infty$, and the difference between $\hat{X}_{n,2}$ and $\hat{X}_{n,2,k}$ is
is asymptotically negligible in probability as $n\ra\infty$ and $k\ra\infty$, i.e.,
\begin{equation}\label{E:Mn2Ce2}
\lim_{n\rightarrow\infty} P\Big( \sup_{0\leq t\leq T, y\geq 0} |\hat{X}_{n,2,k}(t,y)
- \tilde{X}_{n,2,k}(t,y)| > \epsilon \Big) = 0, \quad T >0, \quad \epsilon >0.
\end{equation}
 and
 \begin{equation}\label{E:Mn2Ce1}
\lim\limits_{k\rightarrow\infty} \limsup_{n\rightarrow\infty}P(|\hat{X}_{n,2,k}(t,y) - \hat{X}_{n,2}(t,y)|>\epsilon) = 0, \quad
t ,y \geq  0, \quad \epsilon >0.
\end{equation}

We obtain \eqn{E:Mn2Ce2} easily from Assumption 1 and \eqn{b5} since $\bar{a}$ and $U$ are continuous.
Now we proceed to prove \eqn{E:Mn2Ce1}. We will follow a martingale approach argument similar to the one used in Lemma 5.3 of \cite{KP97},
which relies on their technical Lemma 5.2. Fortunately,  for our two-parameter processes,
the conditions of Lemma 5.2 \cite{KP97} are satisfied by fixing the second argument.
 We have for  $t, y\geq 0$ and $\Upsilon>0$,
\begin{eqnarray}\label{E:Mn2Ce3}
&&P(|\hat{X}_{n,2,k}(t,y) - \hat{X}_{n,2}(t,y)| > \epsilon ) \nonumber\\
&\leq& P(\bar{A}_n(t) > \Upsilon) + P(|\hat{X}_{n,2,k}(t,y) - \hat{X}_{n,2}(t,y)| > \epsilon, \bar{A}_n(t) \leq \Upsilon).
\end{eqnarray}

On $\{\bar{A}_n(t)\leq \Upsilon\}$,
\begin{eqnarray*}
  \hat{X}_{n,2,k}(t,y) - \hat{X}_{n,2}(t,y)
&=&  \int_0^t\int_0^{\infty}(\mathbf{1}^y_{k,t}(s,x)-\mathbf{1}(s+x \leq t +y))d \hat{U}_n(\bar{A}_n(s), F(x)) \nonumber \\
&=&  \frac{1}{\sqrt{n}} \sum_{i = 1}^{A_n(t)\wedge (n \Upsilon)}  \beta_i(\tau_i^n, \eta_i)(t,y),
\end{eqnarray*}
where
\begin{eqnarray*}
 \beta_i(\tau_i^n, \eta_i)(t,y) &=& \sum_{j = 1}^{k} \mathbf{1}(s_{j-1}^k < \tau_i^n \leq s_{j}^k)
\big(\mathbf{1}(t+y - s_j^k < \eta_i \leq t+y - \tau^n_i)  {}\\
&&{}- (F(t+y-\tau_i^n)-F(t+y-s_j^k))  \big) .
\end{eqnarray*}
Define the process $Z_{n}^{(\Upsilon)} \equiv \{Z_{n}^{(\Upsilon)}(t,y): t, y\geq 0\}$ by
\beq
Z_{n}^{(\Upsilon)} (t,y) \equiv \sum_{i = 1}^{A_n(t)\wedge (n\Upsilon)}  \beta_i(\tau_i^n, \eta_i)(t,y), \quad t,y\geq 0.
\eeqno
As in Lemma 5.2  in \cite{KP97}, one can check that for each fixed $y >0$, the process $Z_{n}^{(\Upsilon)}(\cdot, y)
= \{Z_{n}^{(\Upsilon)}(t,y): t\geq 0\}$ is a square integrable martingale with respect to the filtration $\mathbf{F}_{n}
= \{\mathcal{F}_n(t), t\geq 0\}$, where
\begin{eqnarray*}
\mathcal{F}_n(t)&=&  \sigma\{\eta_i \leq s+x: 1\leq i \leq A_n(t), x \geq 0, 0 \leq s \leq t \}  \vee \{A_n(s): 0 \leq s \leq t\}\vee \mathcal{N},
\end{eqnarray*}
and the quadratic variation of  $Z_{n}^{(\Upsilon)}(\cdot, y) $ is
\begin{eqnarray*}
&& \langle Z_{n}^{(\Upsilon)} (\cdot,y)\rangle(t) = \langle Z_{n}^{(\Upsilon)} \rangle(t,y)
=   \sum_{i = 1}^{A^n(t)\wedge (n \Upsilon)} E[ \beta_i(\tau_i^n, \eta_i)(t,y)^2] \nonumber \\
&=& \sum_{i = 1}^{A^n(t)\wedge (n \Upsilon)} \sum_{ j= 1}^k \Big[ \mathbf{1}(s_{j-1} < \tau_i^n \leq s_j^k) (F(t+y-\tau_i^n) -F(t+y-s_j^k)) {}\\
&&{} \cdot(1-(F(t+y-\tau_i^n) - F(t+y-s_j^k)))\Big] \nonumber \\
&\leq &\sum_{i = 1}^{A^n(t)\wedge (n\Upsilon)} \sum_{ j= 1}^k \Big[ \mathbf{1}(s_{j-1}
< \tau_i^n \leq s_j^k) (F(t+y-\tau_i^n) -F(t+y-s_j^k)) \Big] \nonumber \\
&=&\sum_{ j= 1}^k (F(t+y-s_{j-1}^k) -F(t+y-s_j^k))  (A_n(s_j^k) -A_n(s_{j-1}^k)) \nonumber \\
&\leq&  \sup_{1\leq j \leq k} \{ A_n(s_j^k) -A_n(s_{j-1}^k) \},
\end{eqnarray*}
where the last inequality follows from the fact that the sum of the coefficients before the $A_n(s_j^k) -A_n(s_{j-1}^k)$ terms is less than 1.
So for fixed $y \geq 0$,  and on $\{\bar{A}^n(t)\leq \Upsilon\}$,
\begin{eqnarray*}
&&\lim_{k\ra\infty} \limsup_{n\ra\infty}E[ (\hat{X}_{n,2}(t,y) - \hat{X}_{n,2,k}(t,y))^2 ]
=  \lim_{k\ra\infty} \limsup_{n\ra\infty} E\Big[\langle \frac{1}{\sqrt{n}} Z_{n}^{(\Upsilon)}(\cdot,y)\rangle(t)  \Big]  \\
&&\leq \lim_{k\ra\infty} \limsup_{n\ra\infty} E\Big[ \sup_{1\leq j \leq k} \{ \bar{A}_n(s_j^k) -\bar{A}_n(s_{j-1}^k) \} \Big] = 0,
\end{eqnarray*}
where the convergence to 0 holds because of the continuity of $\bar{a}$,
Assumption 1 and $\max_{1\leq j \leq k} (s_j^k - s_{j-1}^k)\rightarrow 0$ as $k\rightarrow \infty$.

Hence, \eqn{E:Mn2Ce3} becomes
\begin{eqnarray*}
P(|\hat{X}_{n,2,k}(t,y) - \hat{X}_{n,2}(t,y)| > \epsilon )
&\leq& P(\bar{A}_n(t) > \Upsilon) + \frac{1}{\epsilon^2} E\Big[\langle \frac{1}{\sqrt{n}} Z_{n}^{(\Upsilon)}(\cdot,y)\rangle(t)  \Big] \\
&\leq& P(\bar{A}_n(t) > \Upsilon) + \frac{1}{\epsilon^2} E\Big[ \sup_{1\leq j \leq k} \{ \bar{A}_n(s_j^k) -\bar{A}_n(s_{j-1}^k) \} \Big].
\end{eqnarray*}

Therefore, by the stochastic boundedness of $\bar{A}_n$, \eqn{E:Mn2Ce1} is proved.
That concludes the demonstration that the finite-dimensional distributions of $(\hat{X}_{n,1}, \hat{X}_{n,2})$ converge to those of $(\hat{X}_1, \hat{X}_2)$ as $n\ra \infty$.~~~\bsq

\section*{Acknowledgments.}   This research was supported by NSF grants DMI-0457095 and CMMI-0948190.

\singlespacing


\begin{thebibliography}{99}




\bibitem{BW72}
Bickel, P.J. and M. J. Wichura.  1971. Convergence criteria for multiparameter stochastic processes
and some applications. {\em Ann. Math. Statist.} 42, 1656--1670.


\bibitem{B68}
Billingsley, P. 1968.
{\em Convergence of Probability Measures}, Wiley
(second edition, 1999).


\bibitem{B67}
Borovkov, A. A. 1967.
On limit laws for service processes in multi-channel systems (in Russian).
{\em Siberian Math J.} 8, 746--763.




\bibitem{B81}
Br\'{e}maud, P. 1981.
{\em Point Processes and Queues:  Martingale Dynamics}, Springer.


\bibitem{Cair}
Cairoli, R. 1972.
 Sur une equation differentielle stochastique.
 {\em Compte Rendus Acad. Sc.} Paris 274 Ser. A, 1739-1742.



\bibitem{CW75}
 Cairoli, R. and J.B. Walsh. 1975.
 Stochastic integrals in the plane.
  { \em Acta Math}. 134, 111-183.



\bibitem{CR81}
Cs\"{o}rg\"{o}, M. and P. R\'{e}v\'{e}sz. 1981.
{\em Strong Approximations in Probability and Statistics},
Wiley, New York.

\bibitem{DM08}
Decreusefond,  L. and P. Moyal. 2008.
A functional central limit theorem for the $M/GI/\infty$ queue.  {\em Ann. Appl. Prob.} Vol. 18, No. 6,  2156-2178.


\bibitem{DW97}
Duffield, N. G. and W. Whitt. 1997.
Control and recovery from rare congestion events in a large multi-server system.
{\em Queueing Systems} 26, 69--104.

\bibitem{EMW93}
Eick, S. G., W. A. Massey and W. Whitt. 1993.
The physics of the $M_t/G/\infty$ queue.
{\em Oper. Res.} 41, 731--742.

\bibitem{EK86}
Ethier, S. N. and T. G. Kurtz.  1986.
{\em Markov Processes: Characterization and Convergence}. Wiley.


\bibitem{GS79}
Gaenssler, P. and W. Stute.  1979.
Empirical processes:  a survey of results for independent and identically distributed random variables.
{\em Ann. Probab.} 7, 193--243.




\bibitem{G82}
Glynn, P. W. 1982.
On the Markov property of the $GI/G/\infty$ Gaussian limit.
{\em Adv. Appl. Prob.} 14, 191--194.

\bibitem{GW91}
Glynn, P. W. and W. Whitt.  1991.
A new view of the heavy-traffic limit theorem for the infinite-server queue.
{\em Adv. Appl. Prob.} 23, 188--209.

\bibitem{GW08}
Goldberg, D. A. and W. Whitt.  2008.
The last departure time from an $M_t/GI/\infty$ queue with a terminating arrival process.
{\em Queueing Systems} 58, 77--104.


\bibitem{HW81}
Halfin, S. and W. Whitt. 1981.
Heavy-traffic limits for queues with many exponential servers.
\emph{Oper. Res.} {\bf 29}(3) 567--587.






\bibitem{HR81}
Harrison, J. M. and M. I. Reiman.  1981.
Reflected Brownian motion in an orthant.
{\em Ann. Probab.} 9, 302--308.




\bibitem{I65}
Iglehart, D. L. 1965.
Limit diffusion approximations for the many server queue and the repairman problem.
{\em J. Appl. Prob.}  2, 429--441.

\bibitem{I73}
Iglehart, D.L. 1973.
Weak convergence of compound stochastic processes.
{\em Stoch. Proc. Appl.} 1, 11--31.

\bibitem{IW70a}
Iglehart, D. L. and W. Whitt. 1970.
Multichannel queues in heavy traffic, I. {\em Adv. Appl. Prob.} 2, 150--177.

\bibitem{IW70b}
Iglehart, D. L. and W. Whitt.  1970.
Multichannel queues in heavy traffic, II:  sequences, networks and batches.
{\em Adv. Appl. Prob.} 2, 355--369..

%

\bibitem{JS87}
Jacod, J. and A. N. Shiryayev.  1987.
{\em Limit Theorems for Stochastic Processes}, Springer.

\bibitem{KR08}
Kang, W. and K. Ramanan. 2010.
Fluid limits of many-server queues with reneging.
To appear in {\em Annals of Applied Probability}. 


\bibitem{KS91}
Karatzas, I. and S. Shreve. 1991. \emph{Brownian Motion and Stochastic
Calculus}. Springer.

\bibitem{KR07}
Kaspi, H. and K. Ramanan. 2010.
Law of large numbers limits for many-server queues.
To appear in {\em Annals of Applied Probability}. 


\bibitem{K02}
Khoshnevisan, D. 2002.
\emph{Multiparameter Processes: An Introduction to Random
Fields}. Springer.



\bibitem{KP97}
Krichagina, E. V. and A. A. Puhalskii.  1997.
A heavy-traffic analysis of a closed queueing system with a $GI/\infty$ service center.
{\em Queueing Systems} 25, 235--280.


\bibitem{KP91}
Kurtz, T.G. and P. Protter.  1991.
Weak limit theorems for stochastic integrals and stochastic differential equations.
{\em Ann. Probab.} 19, 1035--1070.

\bibitem{KP96}
Kurtz, T.G. and P. Protter.  1996.
Weak convergence of stochastic integrals and differential equations II: Infinite dimensional case.
{\em Lecture Notes in Mathematics.}  Vol. 1627, 197-285.


%

\bibitem{L88}
Louchard, G. 1988.
Large finite population queuing systems. Part I:  the infinite server model.
{\em Stochastic Models} 4, 373--505.

\bibitem{M86}
Mamatov, K.M. 1986.
Weak convergence of stochastic integrals with
respect to semimartingales.
{\em Russ. Math. Surv.} 41 (5), 155--156.

\bibitem{MMR98}
Mandelbaum, A., W. A. Massey and M. I. Reiman.  1998.
Strong approximations for Markovian service networks.
{\em Queueing Systems} 30, 149--201.


\bibitem{MW93}
Massey, W. A. and W. Whitt.  1993.
Networks of infinite-server queues with nonstationary Poisson input.
{\em Queueing Systems} 13, 183--250.

\bibitem{MM09}
Mandelbaum, A., P. Momcilovic.
Queues with many servers and impatient customers.
EECS Department, University of Michigan, 2009.





\bibitem{N71}
Neuhaus, G. 1971.
On weak convergence of stochastic processes with multidimensional time
parameter.  {\em Ann. Math. Statist.} 42, 1285--1295.



\bibitem{PTW07}
Pang, G., R. Talreja and W. Whitt. 2007.
Martingale proofs of many-server heavy-traffic limits for Markovian queues.
{\em Probability Surveys}. 4, 193--267.

\bibitem{PW09}
Pang, G., W. Whitt. 2009.
Two-parameter heavy-traffic limits for infinite-server queues: longer version.
Columbia University.  Available at:  http://www.columbia.edu/$\sim$ww2040


\bibitem{PR08}
Puhalskii, A. A. and J.E.  Reed.  2010.
 On many-server queues in heavy traffic.
 {\em Annals of Applied Probability}. Vol. 20, No. 1, 129--195.



\bibitem{R07a}
 Reed, J. E.  2009.
 The G/GI/N queue in the Halfin-Whitt regime I: infinite-server queue system equations.
 {\em Annals of Applied Probability}. Vol. 19, No. 6, 2211--2269.


\bibitem{R07b}
 Reed, J. E. 2007.
The G/GI/N queue in the Halfin-Whitt regime II: idle-time system equations.
working paper, The Stern School, NYU.


\bibitem{RT09}
Reed, J. and R. Talreja. 2009.
Distribution-valued heavy-traffic limits for $GI/GI/\infty$ queues.
{\em Preprint.}


\bibitem{S56}
Skorohod, A. V.  1956.
Limit theorems for stochastic processes.
{\em Prob. Theory Appl.} 1, 261--290.

\bibitem{S71}
 Straf,  M.L. 1971.
Weak convergence of stochastic processes with several
parameters. {\em Proc. Sixth Berkeley Symp. Math. Statist. Prob.}
2, 187--221.

\bibitem{TW08}
Talreja, R. and W. Whitt. 2009.
Heavy-traffic limits for waiting times in many-server queues with abandonments.
 {\em Annals of Applied Probability}. Vol. 19, No. 6, 2137--2175. 

\bibitem{VW96}
van Der Vaart, A. W. and J. Wellner. 1996.
{\em Weak Convergence and Empirical Processes}, Springer.




\bibitem{Walsh}
 Walsh, J.B. 1986.
 Martingales with a multidimensional parameter and stochastic integrals in the plane.
 {\em Lectures in Probability and Statistics}. 329-491, Springer.




\bibitem{W82}
Whitt, W. 1982.
On the heavy-traffic limit theorem for $GI/G/\infty$ queues.
{\em Adv. Appl. Prob.} 14, 171--190.


\bibitem{W02}
Whitt, W. 2002.
{\em Stochastic-Process Limits}. Springer.



\bibitem{WW06}
Whitt, W. 2006.
Fluid models for multiserver queues with abandonments.
{\em Oper. Res.} 54, 37--54.



\bibitem{WW07}
 Whitt, W. 2007.
 Proofs of the martingale FCLT: a review.
 {\em Probability Surveys}. Vol. 4, 268-302.

\bibitem{WZ74}
Wong, E. and M. Zakai, 1974.
Martingales and stochastic integrals for processes with a multidimensional parameter.
{\em Z.Wahrscheinlichkeitstheorie verw. Gebiete}.  29, 109-122.


\bibitem{WZ77}
 Wong, E. and M. Zakai. 1977.
An extension of stochastic integrals in the plane.
{\em Ann. Probability}.  5, 770-778.

\bibitem{Z09}
Zhang, J. 2009.
Fluid models of multi-server queues with abandonment.
{\em Preprint.}



\end{thebibliography}
\end{document}